\documentclass[a4paper]{amsart}
\usepackage{amscd, amssymb}
\usepackage{xypic}
\usepackage[matrix, arrow]{xy}
\usepackage{array}

\def\dec{\operatorname{dec}}
\def\rank{\operatorname{rank}}

\def\A{\mathcal{A}}
\def\B{\mathcal{B}}

\def\E{\mathcal{E}}

\def\hot{\hat{\otimes}}

\newcommand*{\KK}{\mathrm{KK}}

\newcommand{\Z}{\mathbb{Z}}
\newcommand{\C}{\mathbb{C}}

\newcommand{\R}{\mathbb{R}}

\newcommand{\End}{\mathrm{End}}

\def\GL{\operatorname{GL}}
\def\diag{\operatorname{diag}}
\def\Cl{\mathrm{C}l}
\def\Aut{\operatorname{Aut}}

\def\CC{\mathbb C}

\def\E{\mathcal E}

\def\tr{\operatorname{tr}}

\def\lk{\langle}
\def\rk{\rangle}
\def\NN{\mathbb N}
\def\oplus{\bigoplus}
\def\pt{\operatorname{pt}}

\def\RR{\mathbb R}
\def\res{\operatorname{res}}

\def\ZZ{\mathbb Z}

\def\Ad{\operatorname{Ad}}
\def\id{\operatorname{id}}
\def\U{\mathcal U}

\def\eps{\epsilon}

\def\TT{\mathbb T}
\def\dach{{\!\widehat{\ \ }}}

\def\Spin{\operatorname{Spin}}
\def\Pin{\operatorname{Pin}}
\def\O{\operatorname{O}}
\def\SO{\operatorname{SO}}

\def\mod{\operatorname{mod}}

\def\hotimes{{\hot}}

\theoremstyle{plain}
 \newtheorem{theorem}{Theorem}[section]
\newtheorem{corollary}[theorem]{Corollary}
 \newtheorem{lemma}[theorem]{Lemma}

\newtheorem{proposition}[theorem]{Proposition} 
\theoremstyle{definition}

\theoremstyle{definition}
\newtheorem{definition}[theorem]{Definition}

\theoremstyle{remark}
\newtheorem{remark}[theorem]{Remark}

\newtheorem{example}[theorem]{Example}

\numberwithin{equation}{section} \emergencystretch 25pt

\begin{document}

\title[Equivariant $K$-theory of finite dimensional real vector spaces]
  {Equivariant $K$-theory of finite dimensional real vector spaces}

\author[Echterhoff]{Siegfried
  Echterhoff}
   \address{Westf\"alische Wilhelms-Universit\"at M\"unster,
  Mathematisches Institut, Einsteinstr. 62 D-48149 M\"unster, Germany}
\email{echters@uni-muenster.de}

  \author[Pfante]{Oliver
  Pfante}
  \address{Westf\"alische Wilhelms-Universit\"at M\"unster,
  Mathematisches Institut, Einsteinstr. 62 D-48149 M\"unster, Germany}
\email{o\_pfan01@uni-muenster.de}

\begin{abstract} We give a general formula for the equivariant complex $K$-theory $K_G^*(V)$
of a finite dimensional real linear space $V$ equipped with a linear action of a compact group $G$
in terms of the representation theory of a certain double cover of $G$. 
Using this general formula, we give explicit computations in various interesting special cases.
In particular, as an application  we obtain explicit formulas for the $K$-theory of 
$C_r^*(\GL(n,\RR))$, the reduced group C*-algebra of $\GL(n,\RR)$.
\end{abstract}

\maketitle

\section{Introduction}
Let $G$ be a compact group acting linearly on the real vector space $V$. 
In this paper we want to give explicit formulas for the complex equivarant $K$-theory
$K_G^*(V)$ depending on the action of the given group $G$ on $V$.
By use of the positive solution of the Connes-Kasparov conjecture in \cite{CEN},
this will also provide explicit formulas for the $K$-theory
$K_*(C_r^*(H))$ of the reduced group $C^*$-algebra 
$C_r^*(H)$  for any second countable almost connected group $H$
depending on the action of the maximal compact subgroup $G$ of $H$ on 
the tangent space $V=T_{eG}(H/G)$. 

If the action of $G$ on $V$ is orientation preserving (which is always the case if 
$G$ is connected) and lifts to a homomorphism of $G$ to $\Spin(V)$ (or even $\Spin^c(V)$), we get 
the well-known  answer from the equivariant Bott-periodicty theorem.
It implies that 
$$K_G^*(V)\cong K^{*+\dim(V)}_G(\pt)=\left\{\begin{matrix} \oplus_{\rho\in \widehat{G}}\Z & \text{if $*+\dim(V)$ is even}\\
\{0\}& \text{if $*+\dim(V)$ is odd}\end{matrix}\right\}.$$
The obstruction for a linear  action of $G$ on $V$ to lift to a homomorphism into $\Spin(V)$ 
is given by  the Stiefel-Witney class $[\zeta]\in H^2(G,\Z_2)$, where we write $\Z_2:=\Z/2\Z$. 
This class $\zeta$ determines a central group extension
$$
\begin{CD}1@>>> \Z_2@>>> G_\zeta @>>> G @>>> 1.
\end{CD}$$
If we denote by $-1$ the nontrivial element of $\Z_2\subseteq G_\zeta$, then the
irreducible representations of $G_\zeta$ can be devided into the disjoint subsets
$\widehat{G}_\zeta^+$ and $\widehat{G}_\zeta^-$ with
$$\widehat{G}_\zeta^+:=\{\rho\in \widehat{G}: \rho(-1)=1_{V_\rho}\}\quad\quad
\widehat{G}_\zeta^-:=\{\rho\in \widehat{G}: \rho(-1)=-1_{V_\rho}\}.$$
It is then well known (eg., see \cite[\S 7]{CEN}) that $K^*_G(V)$ is a free abelian group 
with one generator for every element $\rho\in \widehat{G}_\zeta^-$.
In particular, it follows that the equivariant 
$K$-theory of $V$ is always concentrated in $\dim(V)$ mod $2$  if the action of $G$ on $V$ is orientation preserving.

The situation becomes more complicated if the action of $G$ on $V$ is not orientation preserving.
As examples show, in this case non-trvial $K$-groups may appear in all dimensions.
The situation has been studied in case of finite groups by Karoubi in \cite{Kar}.
In this paper we use different methods to give a general description of $K_G^*(V)$ 
which works for all compact groups  $G$. We then show in several particular examples how one can extract explicit formulas from our general result, also recovering
the explicit formulas given by Karoubi in the case  of the symmetric group $S_n$ acting
on $\RR^n$ by permuting the coordinates.

To explain our general formula we assume first that $\dim(V)$ is even. Then there is 
 a central extension
$$
\begin{CD}
1 @>>> \Z_2 @>>> \Pin(V) @>\Ad>> \O(V) @>>> 1,
\end{CD}
$$
where $\Pin(V)\subset \Cl_\R(V)$ denotes the Pin-group of $V$ (see \S \ref{sec-prel} for further details
on this extension).
If $\rho:G\to O(V)$ is any continuous homomorphism, let 
$$G_\rho:=\{(x,g)\in \Pin(V)\times G: \Ad(x)=\rho(g)\}.$$
Then $G_\rho$ is a central extension of $G$ by $\ZZ_2$ with quotient map
$q:G_\rho\to G$ given by the projection to the second factor.
%

Let $K_{\rho}\subseteq G_{\rho}$ denote the pre-image of $\SO(V)$ under the homomorphism 
$\rho\circ q:G_\rho\to \O(V)$. Then  $K_\rho=G_\rho$, if the action of $G$ on $V$ is orientation preserving, 
and $[G_\rho:K_\rho]=2$ otherwise. The group $G_\rho$ acts on $\widehat{K}_\rho$ by conjugation.
Let $$\widehat{K}_{\rho}^-=\{\tau\in \widehat{K}_\rho: \tau(-1)=-1_{V_\tau}\}$$
be the set of negative representations of $K_\rho$. This set is invariant under the action 
of $G_\rho$. 
Write $\mathcal O_1$ as the set of all orbits of length one in $\widehat{K}_{\rho}^-$ and 
$\mathcal O_2$ as the set of all  orbits of length two in $\widehat{K}_{\rho}^-$. We then get the 
following general result:

\begin{theorem}\label{thm-K0}
Suppose that $G$ is a compact group acting on the even dimensional real vector space $V$ via the 
homomorphism $\rho: G\to O(V)$. Then, using the above notations, we have:
\begin{enumerate}
\item If $\rho(G)\subseteq \SO(V)$, then 
$$K_0^G(V)\cong \oplus_{[\tau]\in \widehat{G}_\rho^-} \Z \quad\text{and}\quad
K_1^G(V)=\{0\} .$$
\item If $\rho(G)\not\subseteq \SO(V)$, then
$$K_0^G(V)\cong \oplus_{[\tau]\in \mathcal O_2} \Z \quad\text{and}\quad
K_1^G(V)\cong \oplus_{[\tau]\in \mathcal O_1} \Z .$$
\end{enumerate}
\end{theorem}
The odd-dimensional case can easily be reduced to the even case by
passing from $V$ to $V\oplus\R$ together with Bott-periodicity. 
The main tool for proving the theorem is Kasparov's $\KK$-theoretic version of 
equivariant Bott-periodicity, which provides a $\KK^G$-equivalence between $C_0(V)$ and 
the complex Clifford-algebra $\Cl(V)$ (e.g. see \cite[Theorem 7]{Kas3}). By the Green-Julg theorem, this reduces everything to a study of $\KK(\CC, \Cl(V)\rtimes G)$, which then leads to the above representation theoretic description of $K_*^G(V)$.



After having shown the above general theorem we shall consider various special cases in which we 
present more explicit formulas. In particular we shall consider the case of finite groups in 
\S \ref{sec-finite} and 
the case of actions of $\O(n)$ in \S \ref{sec-On} below. In particular, 
we show that for the canonical action of $\O(n)$ on $\RR^n$ we always get $K_{\O(n)}^1(\RR^n)=\{0\}$
and  
$$K_{\O(n)}^0(\RR^n)\cong\left\{\begin{matrix} \ZZ& \text{if $ n=1$}\\
 \oplus_{n\in \NN} \ZZ&\text{if $n>1$}\end{matrix}\right\}$$
 (see Example \ref{ex-cyclic} for the case $n=1$ and Theorem \ref{thm-On} for the case $n>1$).
Another interesting action of $\O(n)$ is the action by conjugation on the space $V_n$ of symmetric matrices in $M(n,\RR)$. By (the solution of) the 
Connes-Kasparov conjecture we have
$K_*\big(C_r^*(\GL(n,\RR))\big)\cong K_{\O(n)}^*(V_n)$ which now allows us to give explicit 
computations of these groups in all cases (see Theorem \ref{thm-GLn} below) 
In particular, in case $n=2$ we get
 $$ K_0\big(C_r^*(\GL(2,\RR))\big)\cong \ZZ\quad\text{and}\quad
 K_1\big(C_r^*(\GL(2,\RR))\big)\cong \oplus_{n\in\NN}\ZZ,$$
which also shows that, for $K_G^*(V)$,  it is possible to have infinite rank in one dimension and finite rank in the other. While these results are stated and prepared in \S \ref{sec-On}, some representation theoretical background and part of the proof are given in  \S \ref{sec-orbits}.

This paper is based on the Diplom thesis of the second named author written under the direction of the first named author at the University of M\"unster. The authors are grateful to Linus Kramer for some useful discussions.

\def\mathcs{{\normalshape\text{C}}^{\displaystyle *}}

\section{Some preliminaries on Clifford algebras}\label{sec-prel}

For the readers convenience, we recall in this section some basic facts on Clifford algebras which will be used throughout this paper. For this let $V$ be a fixed finite dimensional real vector space equipped with an inner product $\lk\cdot,\cdot \rk$. We denote by $\Cl_\R(V)$ the real and by $\Cl(V)=\Cl_\R(V)\otimes_{\R}\C$ the complex Clifford algebras of $V$ with respect to this inner product.
Recall that $\Cl_\R(V)$ is the universal algebra generated by the elements of $V$ subject to the relation
$$v\cdot v=-\lk v,v\rk1.$$
Every element of $\Cl_\R(V)$ is a finite linear combination of elements of the form 
$v_1\cdot v_2\cdots v_k$ with $0\leq k\leq \dim(V)$ and there is a canonical $\ZZ_2$-grading on 
$\Cl_\R(V)$ with grading operator
$$\alpha:\Cl_\R(V)\to \Cl_\R(V); \alpha(v_1\cdots v_k)=(-1)^kv_1\cdots v_k.$$
We shall write $\Cl_\R(V)^0$ and $\Cl_\R(V)^1$ for the even and odd graded elements 
of $\Cl_\R(V)$, respectively.
We also have an involution on $\Cl_\R(V)$ given by $(v_1\cdots v_k)^*=(-1)^kv_k\cdots v_1$. 
With this notation, the Pin-group is defined as
$$\Pin(V)=\{x\in \Cl_\R(V): x^*x=1\;\text{and}\; x v x^*\in V \;\text{for all}\; v\in V\}$$
and  $\Spin(V)=\Pin(V)\cap Cl_\R(V)^0$,
where we regard $V$ as a linear subspace of $\Cl_\R(V)$ in the canonical way.
Similar statements hold for the complex Clifford algebra $\Cl(V)$ if we replace $V$ by its complexification $V_\CC=V\otimes_\R\CC$.  In particular, we obtain the complex Pin-group
$$\Pin^c(V)=\{x\in \Cl(V): x^*x=1\;\text{and}\; x v x^*\in V \;\text{for all}\; v\in V\},$$
and $\Spin^c(V)=\Pin^c(V)\cap \Cl(V)^0$,
where we regard $V$ as a linear subspace of $\Cl(V)$ via the inclusion $x\mapsto x\otimes_\R 1$
of $\Cl_\R(V)$ into $\Cl(V)$. Note that this map also induces an
inclusion  $\iota: \Pin(V)\to \Pin^c(V)$.

If $V$ is even dimensional with dimension $2n$, then $\Cl(V)$ is isomorphic to the full matrix
algebra $M_{2^n}(\CC)$ and the grading on $\Cl(V)$ is given by conjugation
with the element
$$J:=e_1e_2\cdots e_{2n}\in \Spin(V),$$
where $\{e_1,\ldots, e_{2n}\}$ is any given orthonormal base of $V$.
Moreover, there is a short exact sequence
\begin{equation}\label{ex-even}
\begin{CD}
1 @>>> \ZZ_2@>>>\Pin(V) @>\Ad>> \O(V)@>>> 1,
\end{CD}
\end{equation}
where for $x\in \Pin(V)$ the transformation $\Ad(x)\in \O(V)$ is defined by $\Ad(x)(v)=xvx^*$.
The group $\Spin(V)$ is the inverse image of $\SO(V)$ under the adjoint 
homomorphism. For more details on Clifford algebras and the Pin-groups we refer to \cite{ABS}.

Notice that an analogue of the above extension is given in  the odd-dimensional case by 
\begin{equation}\label{ex-odd}
\begin{CD}
1 @>>> \ZZ_2@>>>\Pin(V) @>\A>> \O(V)@>>> 1,
\end{CD}
\end{equation}
where   $\A:\Pin(V)\to \O(V)$ denotes the twisted adjoint given by
$\A(x)(v)=\alpha(x)v x^*$. But this will not play a serious r\^ole in this paper.

Suppose now that  $G$ is a compact group and that $\rho:G\to \O(V)$ is a linear representation of $G$ 
on $V$. Then $\rho$ induces an action $\gamma:G\to \Aut(\Cl(V))$  by defining
$$\gamma_g(v_1\cdots v_k)=\rho(g)(v_1)\cdots \rho(g)(v_k).$$
Note that if 
$x\in \Pin^c(V)$ such that $\Ad(x)=\rho(g)$ for some $g\in G$, then 
\begin{equation}\label{eq-conjugate}
\gamma_g(y)=xyx^*
\end{equation}
for all $y\in \Cl(V)$, which can be checked on the basic elements $x=v_1\cdots v_k$.
Note that this action is compatible with the grading $\alpha=\Ad J$ on $\Cl(V)$.

If $\dim(V)$ is even, we define 
$$G_\rho:=\{(x,g)\in \Pin(V)\times G: \Ad(x)=\rho(g)\}.$$
Then the kernel of the projection $q:G_\rho\to G$  equals $\ZZ_2$ 
and  we obtain 
a central extension
\begin{equation}\label{eq-Grho}
\begin{CD}
1 @>>> \ZZ_2 @>>> G_{\rho} @>q>> G @>>> 1.
\end{CD}
\end{equation}
Let $u:G_\rho\to \Pin(V); u(x,g)= x$ denote the canonical homomorphism. Then it follows from
(\ref{eq-conjugate}) that
\begin{equation}\label{eq-ad}
\gamma_g(y)=u(x,g)yu(x,g)^*\quad\quad\text{for all $y\in \Cl(V)$.}
\end{equation}
Note also that for each $v\in \Pin^c(V)$ we have the equation
$$
JvJ^*=\alpha(v)=\det(\Ad(v))v
$$
which follows from the fact that for $v\in \Pin^c(V)$ the transformation 
$\Ad v$ on $V$ has determinant $1$ if and only if $x\in \Cl(V)^0$.
In particular we get
\begin{equation}\label{eq-grad}
J u(x,g) J^*=\det\rho(g) u(x,g)\quad\quad\text{for all $(x,g)\in G_\rho$}.
\end{equation}
If $\dim(V)$ is odd, we consider the homomorphism $\tilde\rho:G\to \O(V\oplus \R)$
given by $\tilde\rho(g)=\left(\begin{smallmatrix} \rho(g)& 0\\0&1\end{smallmatrix}\right)$.
We then put $G_\rho:=G_{\tilde\rho}$. In this way we obtain a similar central extension as in 
(\ref{eq-Grho}) in the odd-dimensional case.

\section{The main result}\label{sec-main}

Throughout this section we assume that $\rho:G\to \O(V)$ is a linear action of the compact group 
$G$ on the finite-dimensional real vector space $V$ and we let $\gamma:G\to\Aut(\Cl(V))$ denote the 
corresponding action on the complex Clifford algebra $\Cl(V)$.
Let us recall Kasparov's $\KK$-theoretic version of the Bott-periodicity theorem:

\begin{theorem}[{\cite[Theorem 7]{Kas3}}] Let $\rho:G\to \O(V)$ be as above. Then there are classes
$\alpha\in \KK_0^G(C_0(V), \Cl(V))$ and $\beta\in \KK_0^G(\Cl(V), C_0(V))$ which are inverse to each other with respect to the Kasparov product and therefore induce a $\KK^G$-equivalence
between the graded $C^*$-algebra $\Cl(V)$ and the trivially graded algebra $C_0(V)$.
\end{theorem}

From this theorem and the Green-Julg theorem (e.g., see \cite[20.2.7]{Black}), it follows that 
$$K^*_G(V)=\KK_*^G(\C, C_0(V))\cong \KK_*^G(\C, \Cl(V))\cong \KK_*(\C, \Cl(V)\rtimes_\gamma G).$$
So in order to describe the $K$-theory groups  $K^*_G(V)$ it suffices to compute
the groups $\KK_*(\C, \Cl(V)\rtimes_\gamma G)$. Note that we use the $\KK$-notation here and not 
the notation $K_*(\Cl(V)\rtimes_\gamma G)$, since it is important to keep in mind that 
$\Cl(V)\rtimes_\gamma G$ 
is a graded algebra. Indeed, for any function $f\in C(G,\Cl(V))$, regarded as a dense subalgebra 
of $\Cl(V)\rtimes_\gamma G$, the grading operator $\eps:\Cl(V)\rtimes_\gamma G\to \Cl(V)\rtimes_\gamma G$ is given by
$$\eps(f)(g)=\alpha(f(g)),$$
where $\alpha:\Cl(V)\to\Cl(V)$ denotes the grading of $\Cl(V)$.

In case where $\dim(V)=2n$ is even, we shall explicitly describe the crossed product 
$\Cl(V)\rtimes_\gamma G$ as a direct 
sum of full matrix algebras indexed by certain representations of the compact  group 
$G_\rho$ as defined in (\ref{eq-Grho}). 

In general, we have $\Z_2=\{\pm 1\}$ as a central subgroup of $G_\rho$
which  gives us a distinct element $-1\in G_\rho$.
We then write $-g$ for $(-1)g$ 
for all $g\in G_{\rho}$. A function $f\in C(G_\rho)$ is said to be \emph{even} (resp. \emph{odd}), if
$f(-g)=f(g)$ (resp. $f(-g)=-f(g)$) for all $g\in G_\rho$. The even functions can be identified 
with $C(G)$ in a canonical way, and a short computation shows that the convolution on $C(G_\rho)$ restricts to ordinary convolution 
 on $C(G)\subseteq C(G_\rho)$. Similarly, the set of odd functions $C(G_\rho)^-$ is also closed under convolution and 
 involution, and we shall write $C^*(G_\rho)^-$ for its closure in $C^*(G_\rho)$. 
 With this notation we get
 
 \begin{lemma}\label{lem-deco}
The decomposition 
 of $C(G_\rho)$ into even and odd functions induces a direct sum decomposition
  $$C^*(G_\rho)=C^*(G)\oplus C^*(G_\rho)^-$$
with projections 
 $\phi^+:C^*(G_\rho)\to C^*(G), \quad\phi^-:C^*(G_\rho)\to C^*(G_\rho)^-$
 given on $f\in C(G_\rho)$ by
 $$\phi^+(f)(g)=\frac{1}{2}\big(f(g)+f(-g)\big)\quad\text{and}\quad \phi^-(f)(g)=\frac{1}{2}\big(f(g)-f(-g)\big).$$
 There is a corresponding decomposition of $\widehat{G}_\rho$ 
 as a disjoint union $\widehat{G}\cup \widehat{G}_\rho^-$ with 
 $$\widehat{G}_\rho^-:=\{\tau\in \widehat{G}_\rho: \tau(-1)=-1_{V_\tau}\}.$$
 \end{lemma}
 \begin{proof} The proof is fairly straight-forward: If 
 $\tau:G_\rho\to \U(V_\tau)$ is any irreducible unitary representation of 
 $G_\rho$, then $\tau(-1)$ commutes with $\tau(g)$ for all $g\in G_\rho$ and since $\tau(-1)^2=1_{V_\tau}$, it follows from Schur's lemma that either $\tau(-1)=1_{V_\tau}$ or $\tau(-1)=-1_{V_\tau}$.
 In the first case, the representation factors through an irreducible representation of $G$ via
 the quotient map $q:G_\rho\to G$. Thus we obtain a 
 decomposition
 of the dual $\widehat{G}_\rho$ into the disjoint union $\widehat G_\rho=\widehat G\cup \widehat{G}_\rho^-$.
One easily checks that the even (resp. odd) functions on $G_\rho$ are annihilated by (the integrated forms of)
 all elements
of $\widehat{G}_\rho^-$ (resp. $\widehat G$), which then implies that 
this decomposition of $\widehat{G}_\rho$ corresponds to the above described 
decomposition of $C^*(G_\rho)$.
\end{proof} 

\begin{definition}\label{def-negative}
A representation $\tau$ of $G_\rho$ is called \emph{negative}
if $\tau(-1)=-1_{V_\tau}$.
\end{definition}
The negative representations are precisely
those unitary  representation  of $G_\rho$ which factor through 
$C^*(G_\rho)^-$. 
In what follows next, we want to show that in case where $\dim(V)=2n$ is even,
there is a canonical isomorphism
$\Cl(V)\otimes C^*(G_\rho)^-\cong \Cl(V)\rtimes_\gamma G$. To prepare for this result we first note
that the tensor product $\Cl(V)\otimes C^*(G_\rho)^-$ can be realized as the closure of 
the odd functions $f\in C(G_\rho,\Cl(V))$ in the crossed product  $\Cl(V)\rtimes_{\id}G_\rho$
with respect to the trivial action of $G_\rho$ on $\Cl(V)$. The representations are given by
the integrated forms $\varphi\times \sigma$ of
pairs of representations $(\varphi, \sigma)$ on  Hilbert spaces $H$ such that
$\varphi$ is a $*$-representation of $\Cl(V)$, $\sigma$ is a negative 
unitary representation of $G_\rho$, and
$$ \varphi(x)\sigma_g=\sigma_g\varphi(x)\quad\forall g\in G_\rho, x\in \Cl(V).$$
For the following proposition recall that $q:G_\rho\to G$ denotes the quotient map.

\begin{proposition}\label{prop-iso}
Suppose that $\dim(V)=2n$ is even. Then there is a canonical isomorphism 
$$\Theta: \Cl(V)\otimes C^*(G_\rho)^-\stackrel{\cong}{\longrightarrow} \Cl(V)\rtimes_\gamma G$$
which sends an odd function $f\in C(G_\rho,\Cl(V))$ to the function $\Theta(f)\in C(G,\Cl(V))\subseteq \Cl(V)\rtimes_\gamma G$,
given 
by
$$\Theta(f)(q(g))=\frac{1}{2}\big(f(g)u_g+f(-g)u_{-g}\big).$$
Under this isomorphism, a representation $\varphi\times \sigma$ of
$\Cl(V)\otimes C^*(G_\rho)^-$ on a Hilbert space $H$ 
corresponds to the representation $\varphi\times \tau$ 
of $ \Cl(V)\rtimes_\gamma G$ on $H$ with $\tau:G\to \U(H)$ given by
$$\tau(q(g))=\varphi(u_g^*)\sigma(g).$$
\end{proposition}
\begin{proof}
 For the proof
we first inflate the action   $\gamma: G\to \Aut(\Cl(V))$ to an action 
$\tilde\gamma:G_\rho\to \Aut(\Cl(V))$ in the obvious way. It follows then from 
(\ref{eq-ad}) that this action
 is implemented by the canonical homomorphism
$u:G_\rho\to \Pin(V)$
in such a way that 
$$\tilde\gamma_g=\Ad u_g$$
for all $g\in G_\rho$. It follows that the crossed product 
$\Cl(V)\rtimes_{\tilde\gamma} G_\rho$ is isomorphic to  
$\Cl(V)\rtimes_{\id} G_\rho\cong C^*(G_\rho)\otimes\Cl(V)$
with isomorphism $\Phi:\Cl(V)\rtimes_{\id} G_\rho \to \Cl(V)\rtimes_{\tilde\gamma} G_\rho$
given on the dense subalgebra $C(G_\rho, \Cl(V))$ by
$$\big(\Phi(f)\big)(g)=f(g)u_g.$$
On the other hand, we have a canonical surjective $*$-homomorphism
 $$\Psi: \Cl(V)\rtimes_{\tilde\gamma}G_\rho\to \Cl(V)\rtimes_\gamma G$$
 given on $C(G_\rho,\Cl(V))$ by $\Psi(f)(q(g))=\frac{1}{2}(f(g)+f(-g))$.
 
 We claim that the $*$-homomorphism $\Theta: \Cl(V)\otimes C^*(G_\rho)\to \Cl(V)\rtimes_\gamma G$ given by the composition
 $$\begin{CD}
\Cl(V) \otimes  C^*(G_\rho) \cong\Cl(V)\rtimes_{\id} G_\rho
 @>\Phi>> \Cl(V)\rtimes_{\tilde\gamma}G_\rho
 @>\Psi>>\Cl(V)\rtimes_\gamma G
 \end{CD}
 $$
 factors through the desired isomorphism $\Cl(V) \otimes C^*(G_\rho)^- \cong \Cl(V)\rtimes_{\gamma}G$.
 It is then clear that it is given by the formula as in the proposition. 

For the proof of the claim it suffices to show that the map
 $$\widehat{\Theta}:(\Cl(V)\rtimes_\gamma G)\dach\to (\Cl(V)\otimes C^*(G_\rho))\dach; 
 \varphi\times \tau\mapsto  (\varphi\times \tau)\circ \Theta$$
 is injective with image $(\Cl(V)\otimes C^*(G_\rho)^-)\dach$.  

To see this let  $(\varphi, \tau)$ be any covariant representation of $\big(\Cl(V),G,\gamma\big)$ 
on the Hilbert space $H$.
Composing it with the quotient map $\Psi$ gives the covariant  representation 
$(\varphi, \tau\circ q)$ of $\big(\Cl(V), G_\rho,\tilde\gamma\big)$.
This representation corresponds to the representation $(\varphi,\sigma)$ of $\big(\Cl(V), G_\rho,\id\big)$
with $\sigma:G_\rho\to \U(V_\tau)$ given by
$$\sigma_g=\varphi(u_g)\cdot \tau_{q(g)}.$$
To see this, we 
simply compute for given $f\in C(G_\rho, \Cl(V))$ the integrated form
\begin{align*}
\varphi\times\tau\circ q(\Phi(f))&=\int_{G_\rho} \varphi(\Phi(f)(g))\tau_{q(g)}\,dg\\
&=\int_{G_\rho} \varphi(f(g)u_g)\tau_{q(g)}\,dg=
\varphi\times \sigma(f).
\end{align*}
Since $u_{-1}=-1$ in $\Cl(V)$, we see that $\sigma(-1)=\varphi(-1)\tau(1)=-1_{H}$, which implies that 
$\sigma$ 
factors through a representation of $C^*(G_\rho)^-$, and therefore $\varphi\times\sigma=(\varphi\times \tau)\circ \Theta$
is an irreducibe representation of $\Cl(V)\otimes C^*(G_\rho)^-$. 
Conversely, assume that a representation  $\varphi\times \sigma$ of $\Cl(V)\otimes C^*(G_\rho)^-$ 
on a Hilbert space $H$ is given.
Then one checks that $(\varphi,\tau)$ with $\tau(q(g))=\varphi(u_g^*)\sigma(g)$ is a 
covariant representation of $\big(\Cl(V), G,\gamma\big)$ such that 
$\varphi\times\sigma=(\varphi\times \tau)\circ \Theta$, and the result follows.
\end{proof}

 Recall that by the Peter-Weyl Theorem the C*-algebra $C^*(G_\rho)$ of the compact group $G_\rho$
 has a decomposition 
 $$C^*(G_\rho)=\oplus_{\tau\in \widehat{G}_\rho}\End(V_\tau)$$
 with projection from $C^*(G_\rho)$ onto the summand $\End(V_\tau)$ given by 
 $f\mapsto \tau(f)=\int_{G_\rho} f(g)\tau_g\,dg$.
 The above decomposition 
  together with Proposition \ref{prop-iso} induces a decomposition
 $$\Cl(V)\rtimes_{\gamma}G\cong \oplus_{\tau\in \widehat{G}_\rho^-} \Cl(V)\hotimes \End(V_\tau).$$
 We need to analyze the grading on the direct sum decomposition  induced  by the grading $\eps$ 
 of $\Cl(V)\rtimes_\gamma G$. Recall that the latter is given on functions $f\in C(G,\Cl(V))$ by
$$\eps(f)(g)=\alpha(f(g))=Jf(g)J^*$$
with $J=e_1\cdots e_{2n}$ and $\{e_1,\ldots, e_{2n}\}$ an orthonormal basis of $V$.
If $\Theta$ is the isomorphism of 
Proposition \ref{prop-iso}
we compute for any elementary tensor $x\otimes f\in \Cl(V)\otimes C(G_\rho)^-\subset C(G_\rho,\Cl(V))$
\begin{align*}
\eps(\Theta(x\otimes f))(q(g))&=J\Theta(x\otimes f)(q(g) )J^*\\
&=J \frac{1}{2}\big(xf(g)u_g+xf(-g)u_{-g}\big) J^*\\
&=\frac{1}{2}\big(JxJ^*f(g)Ju_gJ^*+JxJ^*f(-g)Ju_{-g}J^*\big) \\
&=JxJ^*\det(\rho(q(g)))\frac{1}{2}\big(f(g)u_g+f(-g)u_{-g}\big)\\
&=\Theta(\alpha(x)\otimes (\det\circ \rho)\cdot f)(q(g)),
\end{align*}
where  the second to last equation follows from (\ref{eq-grad}). This shows that the grading 
on $\Cl(V)\hotimes C^*(G_\rho)^-$ corresponding to $\eps$ is the diagonal grading 
given by the standard grading $\alpha$ on $\Cl(V)$ and the grading 
on $C^*(G_\rho)^-$ given on $C(G_\rho)^-$ via point-wise multiplication with the 
$\ZZ_2$-valued character 
$$\mu:G_\rho\to \ZZ_2; \mu(g)=\det\circ \rho\circ q(g).$$
Now, for any function $f\in C(G_\rho)$ and $\tau\in \widehat{G}_\rho$ we get
$$
\tau(\mu\cdot f)(g)=\int_{G_\rho} \mu(g) f(g) \tau_g\,dg=(\mu\cdot\tau)(f),
$$
which implies that the corresponding grading on 
$C^*(G_\rho)^-=\oplus_{\tau\in \widehat{G}_\rho^-}\End(V_\tau)$
induces an inner automorphism on the block $\End(V_\tau)$ if $\tau\cong \mu\cdot\tau$
and  intertwines $\End(V_\tau)$ with $\End(V_{\mu\tau})$ if $\mu\tau\not\cong  \tau$. 

Write $M_\tau:=\Cl(V)\otimes \End(V_\tau)$. Then $M_\tau$ is isomorphic to a full matrix
algebra. If $\tau\cong \mu\cdot\tau$, this summand of $\Cl(V)\otimes C^*(G_\rho)^-$ 
is fixed by the grading  and $M_\tau$ is Morita equivalent (as graded algebra) to the trivially
graded algebra $M_\tau$. 
If $\mu\tau\not\cong  \tau$, the grading intertwines $M_\tau$ with $M_{\mu\tau}$ and the direct
sum $M_\tau\oplus M_{\mu\tau}$ is isomorphic to the algebra $M_\tau\oplus M_\tau$
with the standard odd grading given by $(S,T)\mapsto (T,S)$. Thus 
$M_\tau\oplus M_{\mu\tau}$ is isomorphic to $M_\tau\hotimes\Cl_1$, where 
$\Cl_1=\C\oplus\C$ denotes the first Clifford algebra. We therefore obtain a decomposition
\begin{equation}\label{eq-deco}
\Cl(V)\rtimes_\gamma G\cong 
\left(\oplus_{\tau\in \mathcal O_1} M_\tau\right)
\oplus \left(\oplus_{\{\tau,\mu\tau\}\in \mathcal O_2} M_\tau\hotimes \Cl_1\right),
\end{equation}
where $\mathcal O_1$ denotes the set of fixed points in $\widehat{G}_\rho^-$ under the
order-two transformation $\tau\mapsto \mu\tau$, 
and $\mathcal O_2$ denotes the set of orbits of length two under this action.
Using this, it is now easy to prove:

\begin{theorem}\label{thm-K}
Let $\rho:G\to \O(V)$ be a linear action of the compact group $G$ on the finite-dimensional
real vector space $V$ and let $\gamma:G\to\Aut(\Cl(V))$ denote the corresponding 
action on $\Cl(V)$.
Let   $\widehat{G}_\rho^-$, $\mathcal O_1$ and $\mathcal O_2$ be as above.
 Then 
$$K_G^*(V)\cong \KK_*(\C, \Cl(V)\rtimes_\gamma G)= \left\{\begin{matrix}\oplus_{\tau\in \mathcal O_1} \ZZ&\text{if $*+\dim(V)=0\mod 2$}\\
\oplus_{\{\tau,\mu\tau\}\in \mathcal O_2} \ZZ&\text{if $*+\dim(V)=1\mod 2$}\end{matrix}\right\}.$$
\end{theorem}
\begin{proof} We first assume that $\dim(V)=2n$ is even.
The isomorphism $K_G^*(V)\cong \KK_*(\C, \Cl(V)\rtimes_\gamma G)$ is Kasparov's
Bott-periodicity theorem. The decomposition of $\Cl(V)\rtimes_\gamma G$ of (\ref{eq-deco}) implies a decomposition
\begin{align*}
\KK_*(\C, \Cl(V)\rtimes_\gamma G)&\cong \KK_*\big(\C,\oplus_{\tau\in \mathcal O_1}M_\tau\big)\oplus
\KK_*\big(\C,\oplus_{\{\tau,\mu\tau\}\in \mathcal O_2}M_\tau\hotimes\Cl_1\big)\\
&=K_*\big(\oplus_{\tau\in \mathcal O_1}M_\tau\big)\oplus K_{*+1}\big(\oplus_{\{\tau,\mu\tau\}\in \mathcal O_2}M_\tau\big)\\
&=\left(\oplus_{\tau\in \mathcal O_1}K_*(M_\tau)\right)\oplus\left( \oplus_{\{\tau,\mu\tau\}\in \mathcal O_2}K_{*+1}(M_\tau)\right),
\end{align*}
and the result follows from the fact that $K_0(M)=\ZZ$ and $K_1(M)=\{0\}$  for any full matrix algebra 
$M$.

If $\dim(V)$ is odd,  we defined $G_\rho=G_{\tilde\rho}$ where
$\tilde\rho:G\to\O(V\oplus \R)$, $\tilde\rho(g)=\left(\begin{smallmatrix} \rho(g)& 0\\ 0&1\end{smallmatrix}
\right)$. Note that we have $K_G^*(V)\cong K_G^{*+1}(V\oplus \R)$ by Bott-periodicity.
Thus the 
result follows from the even-dimensional case applied to the action $\tilde\rho$ on $V\oplus\R$.
\end{proof}

Note that the homomorphism $\mu$ is trivial if and only if  $\rho$ takes image in $\SO(V)$, i.e., the action of
$G$ on $V$ is orientation preserving.  In this case we get $\mathcal O_1=\widehat{G}_\rho^-$ and 
$\mathcal O_2=\emptyset$. Moreover, we then have
$$G_\rho=\{(x,g)\in \Spin(V)\times G: \Ad x=\rho(g)\}.$$  
This construction of $G_\rho$ makes also sense if the dimension of $V$ is odd,
and it then coincides with the extension $G_{\tilde\rho}$ we obtain by passing to 
the action $\tilde\rho$ of $G$ on $V\oplus \R$. We leave the verification of this simple
fact to the reader.
We therefore recover the following well-known result (e.g., see \cite[\S7]{CEN}): 

\begin{corollary}\label{cor-oriented}
Let $\rho:G\to \SO(V)$ be an orientation preserving 
linear action of the compact group $G$ on the finite dimensional
real vector space $V$. Then 
$$K_G^*(V)\cong  \left\{\begin{matrix}\oplus_{\tau\in  \widehat{G}_\rho^-} \ZZ&\text{if $*+\dim(V)=0 \mod 2$}\\
\{0\} &\text{if $*+\dim(V)=1\mod 2$}\end{matrix}\right\}.$$
\end{corollary}

At this point it might be interesting to notice, that for non-finite compact groups $G$ 
the cardinality of $\widehat{G}_\rho^-$ is always countably infinite if $G$ is second countable.
This follows from

\begin{lemma}\label{lem-reps} 
Let $1\to \ZZ_2\to H \to G\to 1$ be a central extension of an infinite compact second countable 
group $G$ by $\ZZ_2$. Then the set $\widehat{H}^-$ of equivalence classes of negative irreducible representations of $H$ is countably infinite.
\end{lemma}
\begin{proof}
We decompose $L^2(H)$ as a direct
sum $L^2(G)\oplus L^2(H)^-$, where we identify $L^2(G)$ with the set of even functions 
in $L^2(H)$, and where $L^2(H)^-$ denotes the set odd functions in $L^2(H)$.
Since $G$ is not finite, both spaces are separable infinite dimensional Hilbert spaces.
The regular representation $\lambda_H:H\to \U(L^2(H))$ then decomposes 
into the direct sum $\lambda_G\oplus\lambda_H^-$, with $\lambda_H^-$ a negative 
representation of $H$.  By the Peter-Weyl Theorem 
we get a decomposition
$$\lambda_G{\oplus}\lambda_H^-=\oplus_{\tau\in \widehat{H}}d_\tau\cdot\tau,$$
where $d_\tau$ denotes the dimension of $\tau$ and $d_\tau\cdot \tau$ stands for the 
$d_\tau$-fold direct sum of $\tau$ with itself.
Since a direct summand of a negative (resp. positive)  representation must be negative (resp. postive), we get the decomposition 
$$\lambda_H^-=\oplus_{\tau\in \widehat{H}^-}d_\tau\tau.$$
The result now follows from the fact that all representations $\tau$ in this decomposition are 
finite dimensional.
\end{proof}

We proceed with a discussion of the special case, where the 
homomorphism $\rho:G\to \O(V)$ (resp. $\tilde\rho:G\to \O(V\oplus \R)$ in case where 
$\dim(V)$ is odd) factors through a homomorphism $v:G\to \Pin^c(V)$ (resp. $\Pin^c(V\oplus\R)$).
We then say that the action satisfies a \emph{$\Pin^c$-condition}.
If  $\rho(G)\subseteq \SO(V)$, this implies that the action factors through a homomorphism 
to $\Spin^c(V)$, in which case we say that the action staisfies a \emph{$\Spin^c$-condition}.
Assume first that $\dim(V)$ is even. It follows then from (\ref{eq-conjugate}) that 
the corresponding action on $\Cl(V)$  is given by $\gamma=\Ad v$,
which then implies that
$$\Cl(V)\rtimes_{\gamma}G\cong \Cl(V)\rtimes_{\id}G\cong \Cl(V)\otimes C^*(G)$$
with isomorphism $\Theta: \Cl(V)\rtimes_{\id}G\to \Cl(V)\rtimes_{\gamma}G$ given on functions 
$f\in C(G,\Cl(V))$ by 
$$\big(\Theta(f)\big)(g)=f(g)v_g.$$
Let $\eps$ be the grading on $\Cl(V)\rtimes_\gamma G$.
On elementary tensors $x\otimes f\in \Cl(V)\otimes C(G)\subset \Cl(V)\otimes C^*(G)$
we compute
\begin{align*}
 \eps\big(\Theta(x\otimes f)\big)(g)&=J(xf(g)v_g)J^*=JxJ^* \det(\rho(g))f(g)v_g\\
 &=
\Theta\big(\alpha(x)\otimes (\det\circ \rho) \cdot f\big)(g),
\end{align*}
where the second to last equation follows from (\ref{eq-grad}). So we see that if $\mu:G\to \ZZ_2$ 
denotes the character $\mu=\det\circ \rho$ then the grading on $\Cl(V)\rtimes_\gamma G$ 
corresponds to the diagonal grading on $\Cl(V)\otimes C^*(G)$ given on the second factor
by multiplication with the character $\mu$. Passing to $V\oplus \R$ in case where 
$\dim(V)$ is odd, we now obtain the following theorem, where the proof proceeds
precisely as in Theorem \ref{thm-K}:

\begin{theorem}\label{thm-pinc}
Assume that the linear action $\rho:G\to \O(V)$ satisfies a $\Pin^c$-condition as defined above.
Let $\mathcal O_1$ and $\mathcal O_2$ denote the sets of orbits in $\widehat{G}$
under the order two  transformation $\tau\mapsto \mu\tau$ with $\mu=\det\circ\rho:G\to \ZZ_2$.
 Then 
$$K_G^*(V)\cong \KK_*(\C, \Cl(V)\rtimes_\gamma G)= \left\{\begin{matrix}\oplus_{\tau\in \mathcal O_1} \ZZ&\text{if $*+\dim(V)=0\mod 2$}\\
\oplus_{\{\tau,\mu\tau\}\in \mathcal O_2} \ZZ&\text{if $*+\dim(V)=1\mod 2$}\end{matrix}\right\}.$$
In particular, if $\rho$ satisfies a $\Spin^c$-condition, then $\mathcal O_1=\widehat{G}$ and $\mathcal O_2=\emptyset$ and then
$$K_G^*(V)\cong \KK_*(\C, \Cl(V)\rtimes_\gamma G)= \left\{\begin{matrix}\oplus_{\tau\in \widehat{G}} \ZZ&\text{if $*+\dim(V)=0\mod 2$}\\
\{0\}&\text{if $*+\dim(V)=1\mod 2$}\end{matrix}\right\}.$$
\end{theorem}

Note that the above statements are different from the statement given in the introduction, where
we describe the $K$-groups in terms of the kernel $K_\rho$ of the character $\mu:G_\rho\to \Z_2$ 
and the action of $G_\rho$ on $\widehat{K}_\rho$. 

In order to  see that the above results can be formulated as in the introduction let us assume 
that $G$ is any compact group, $\mu:G\to \Z_2$ is a non-trivial 
continuous group homomorphism, and $K:=\ker\mu\subseteq G$. Then 
$K$ is a normal subgroup of index two in $G$ and $G$ acts on $\widehat{K}$ by conjugation.
This action is trivial on $K$ and therefore factors to an action of $G/K\cong \Z_2$ on $\widehat{K}$.
Thus the $G$-orbits in $\widehat{K}$ are either of length one or of length two. 

\begin{proposition}\label{prop-restrict}
Let $\mu:G\to \Z_2$ and $K=\ker\mu$ as above. 
Consider the action of $\Z_2$ on $\widehat{G}$ given on the generator by multiplying with 
$\mu$ and let $\Z_2\cong G/K$ act on $\widehat{K}$ via conjugation.
Then there is a canonical bijection
$$\res:\widehat{G}/\Z_2\to \widehat{K}/\Z_2$$
 which maps orbits of length one in $\widehat{G}$ to 
orbits of length two in $\widehat{K}$ and vice versa. 
\end{proposition}

This proposition is well known to the experts and follows from basic representation theory.
An elementary proof 
in case of finite groups $G$ is given in \cite[Theorem 4.2 and Korollar 4.3]{HH} and the same
arguments work for general compact groups. Let $g\in G\smallsetminus K$ be any fixed
element. The basic steps for the proof are as follows:
\begin{itemize}
\item If $\tau$ is an irreducible representation of $G$, the restriction $\tau|_K$ is either
irreducible, or it decomposes into the direct sum
 $\sigma\oplus g\cdot \sigma$ for some $\sigma\in \widehat{K}$ with $g\cdot\sigma(k)=\sigma(g^{-1}kg)$
 for $k\in K$.
 \item $\tau|_K$ is irreducible if and only if $\tau\not\cong \mu\tau$.
 \item If $\tau|_K\cong \sigma\oplus g\cdot \sigma$ for some $\sigma\in \widehat{K}$,
 then $\sigma\not \cong g\cdot \sigma$.
 \end{itemize}
The map $\res:\widehat{G}/\Z_2\to \widehat{K}/\Z_2$ of the proposition is then given by
sending the an orbit $\{\tau,\mu\tau\}$ of length two in $\widehat{G}$ to the 
orbit $\{\tau|_K\}$ of length one in $\widehat{K}$ and an orbit $\{\tau\}$ of length one
in $\widehat{G}$ to the orbit $\{\sigma, g\cdot\sigma\}$ of length two in $\widehat{K}$ 
determined by $\tau|_K\cong \sigma\oplus g\cdot\sigma$. It follows  from Frobenius reciprocity that
this map is onto.

We now come back to the description of $K_G^*(V)$:

\begin{corollary}\label{cor-K}
Suppose that $G$ is a compact group and $\rho:G\to\O(V)$ is a linear action of $G$
on the finite dimensional real vector space $V$. Assume that
$\mu=\det\circ \rho\circ q:G_\rho\to \Z_2$ is not trivial, i.e the action of $G$ on $V$ is not orientation preserving. Let $K_\rho=\ker\mu\subseteq G$ and let $$\widehat{K}_\rho^-=\{\sigma\in \widehat{K}_\rho: \sigma(-1)=- 1_{V_\sigma}\}.$$
Then $\widehat{K}_\rho^-$ is invariant under the conjugation action of 
$G_\rho$ on $\widehat{K}_\rho$. 
Let $\mathcal O_1$ and $\mathcal O_2$ denote the set of orbits  of length one or two in
$\widehat{K}_\rho^-$. Then 
$$K_G^*(V)\cong \left\{\begin{matrix}\oplus_{\sigma\in \mathcal O_1} \ZZ&\text{if $*+\dim(V)=1\mod 2$}\\
\oplus_{\{\sigma,g\sigma\}\in \mathcal O_2} \ZZ&\text{if $*+\dim(V)=0\mod 2$}\end{matrix}\right\}.$$
\end{corollary}
\begin{proof} We first note that the central subgroup $\Z_2$ of $G_\rho$ lies in the kernel 
of $\mu$, since $\mu$ factors through $G$. So the definition of 
$\widehat{K}_\rho^-$ makes sense. Moreover, since $\Z_2$ is central in $G_\rho$, we get 
$g(-1)g^{-1}=-1$ for all $g\in G_{\rho}$, which implies that $\widehat{K}_\rho^-$ 
is invariant under the conjugation action. The description of the bijection 
$\res:\widehat{G}_\rho/\Z_2\to \widehat{K}_\rho/\Z_2$ of Proposition \ref{prop-restrict} 
given above now implies that it restricts to a bijection 
$\res^-:\widehat{G}_\rho^-/\Z_2\to \widehat{K}_\rho^-/\Z_2$, which maps 
orbits of length one to orbits of length two and vice versa. The result  then follows 
directly from Theorem \ref{thm-K}.
\end{proof}

In case where $\rho:G\to\O(V)$ satisfies a $\Pin^c$-condition as considered in Theorem \ref{thm-pinc},
the same proof as for the above corollary together with Theorem \ref{thm-pinc} gives

\begin{corollary}\label{cor-Kpinc}
Suppose that the linear action $\rho:G\to\O(V)$ satisfies a $\Pin^c$-condition and
assume that $\mu=\det\circ \rho:G\to \ZZ_2$
is non-trivial.  Let $K=\ker\mu\subseteq G$ and let
 $\mathcal O_1$ and $\mathcal O_2$ denote the set of $G$-orbits in $\widehat{K}$ 
of length one and two, resectively. Then 
$$K_G^*(V)\cong \left\{\begin{matrix}\oplus_{\sigma\in \mathcal O_1} \ZZ&\text{if $*+\dim(V)=1\mod 2$}\\
\oplus_{\{\sigma,g\sigma\}\in \mathcal O_2} \ZZ&\text{if $*+\dim(V)=0\mod 2$}\end{matrix}\right\}.$$
\end{corollary}

\section{Actions of finite groups}\label{sec-finite}

In this section we want to study the case of finite groups in more detail. This case was already considered by Karoubi in \cite{Kar}, but the methods used 
here are different from those used by Karoubi.
We first notice that it follows from Theorem \ref{thm-K} that for actions of finite groups $G$, the $K$-theory groups $K_G^*(V)$ are always 
finitely generated free abelian groups. In what follows next we want to give formulas 
for the ranks of these groups in terms of conjugacy classes. 
Recall that for every finite group $G$ the number $|\widehat{G}|$ of
(equivalence classes of) irreducible representations of $G$ equals the number $C_G$ of conjugacy classes in $G$. 

We first look at the case where the action $\rho:G\to \O(V)$ 
satisfies a $\Pin^c$-condition.
Let $\mu=\det\circ \rho:V\to \ZZ_2$. If $\mu$ is trivial
it follows from Theorem \ref{thm-pinc} that $\rank(K_G^*(V))=|\widehat{G}|=C_G$ if 
$*+\dim(V)=0 \mod 2$ and $\rank(K_G^*(V))=0$ else.

If $\mu$ is non-trivial, let $K=\ker\mu$ and let
$O_1$ and $O_2$ denote the number of $G$-orbits in $\widehat{K}$ of length one and two, respectively. 
The numbers $O_1,O_2$ satisfy
the equations
\begin{equation}\label{eq-O1}
O_1+2O_2=|\widehat{K}|=C_K\quad\text{and}\quad 2O_1+O_2=|\widehat{G}|=C_G.
\end{equation}
Indeed, the first equation follows from the obvious fact that $O_1+2O_2$ coincides with the
number of irreducible representations of $K$ and it follows from Proposition \ref{prop-restrict}
that $2O_1+O_2$ coincides with the number of irreducible representations of ${G}$.
By basic linear algebra the equations (\ref{eq-O}) have the unique solutions 
\begin{equation}\label{eq-sol}
O_1=\frac{1}{3}\big(2C_G-C_K)\quad\text{and}\quad O_2=\frac{1}{3}(2C_K-C_G).
\end{equation}
Combining all this with Corollary \ref{cor-Kpinc} implies

\begin{proposition}\label{prop-finitepinc}
Suppose that the linear action $\rho:G\to \O(V)$ of the finite group $G$
 satisfies a $\Pin^c$-condition.  Then 
 \begin{equation}\label{eq-spin}
\rank(K_G^*(V))=\left\{\begin{matrix} C_G& \text{if $*+\dim(V)=0 \mod 2$}\\
0 &\text{if $*+\dim(V)=1 \mod 2$}\end{matrix}\right\}
\end{equation}
if the action is orientation preserving, and 
 \begin{equation}\label{eq-Opin}
\rank(K_G^*(V))=\left\{\begin{matrix}\frac{1}{3}\big(2C_K-C_G) & \text{if $*+\dim(V)=0 \mod 2$}\\
\frac{1}{3}(2C_G-C_K) &\text{if $*+\dim(V)=1 \mod 2$}\end{matrix}\right\}
\end{equation}
if the action is not orientation preserving, where $K=\{g\in G: \det(\rho(g))=1\}$.
\end{proposition}

\begin{example}\label{ex-cyclic}
Let $G=\ZZ_m$, the cyclic group of order $m$, and let $\rho:\ZZ_m\to \O(V)$ be a linear action 
of $\ZZ_m$ on $V$. We claim that $\rho$ automatically satisfies a $\Pin^c$-condition.
For this let  $g$ be a generator of $\ZZ_m$. By passing to $V\oplus \R$ if necessary, we may assume without
loss of generality that $\dim(V)$ is even. Choose $u\in \Pin^c(V)$ such that $\Ad(c)=\rho(g)$.
Then $\Ad(u^m)=\rho(g^m)=\rho(e)=1_V$, and there exists $\lambda \in \TT$ with $u^m=\lambda 1$.
Changing $u$ into $\zeta u$, where $\zeta\in \TT$ is an $m$th root of $\bar{\lambda}$, we 
obtain a well defined homomorphism $v:\ZZ_m\to \Pin^c(V)$ which sends $g$ to $\zeta u$ 
such that $\Ad v(g^k)=\rho(g^k)$ for all $k\in \ZZ$.  If the action $\rho$ takes image in $\SO(V)$
(which is automatic if $m$ is odd), we get 
$$K_{\ZZ_m}^*(V)\cong \left\{\begin{matrix}\ZZ^m&\text{if $*+\dim(V)=0\mod 2$}\\
\{0\} &\text{if $*+\dim(V)=1\mod 2$}\end{matrix}\right\},$$
since $C_G=m$. Assume now that $m$ is even and the action is not orientation preserving.
Then the group $K=\ker(\det\circ \rho)$ is cyclic
of order $\frac{m}{2}$.  With the values 
$C_G=m$ and $C_K=\frac{m}{2}$, Proposition \ref{prop-finitepinc} implies
$$K_{\ZZ_m}^*(V)\cong
 \left\{\begin{matrix} \ZZ^{m/2}&\text{if $*+\dim(V)=1\mod 2$}\\
\{0\}&\text{if $*+\dim(V)=0\mod 2$}\end{matrix}\right\}.$$
As a particular example, consider the action of $\ZZ_2$ on $\RR$ by reflection. Then we get
$$K_{\ZZ_2}^0(\RR)=\ZZ\quad\text{and}\quad K_{\ZZ_2}^1(\RR)=\{0\}.$$
\end{example}
\medskip

We now study the general case. Consider the 
central extension $1\to\ZZ_2\to G_\rho\stackrel{q}{\to}G\to 1$ 
as in (\ref{eq-Grho}). As shown in the previous section we have a disjoint decomposition 
$\widehat{G}_\rho=\widehat{G}_\rho^-\cup\widehat{G}$.  We therefore get 
$$|\widehat{G}_\rho^-|=|\widehat{G}_\rho|-|\widehat{G}|=C_{G_\rho}-C_G.$$
Thus, if the action $\rho:G\to \SO(V)$ is orientation preserving,  it follows from 
Corollary \ref{cor-oriented} that 
$\rank(K_G^*(V))=C_{G_\rho}-C_G$ if 
$*+\dim(V)=0 \mod 2$ and $\rank(K_G^*(V))=0$ else.

If $\rho$ is not orientation preserving we obtain the non-trivial character
$\mu :G_\rho\to \ZZ_2$, $\mu=\det\circ \rho\circ q$. Let $K_\rho=\ker\mu$ and let 
$O_1$ and $O_2$ denote the number of $G_\rho$-orbits of length one and two in $\widehat{K}_\rho^-$, respectively.
Similar as for $|\widehat{G}_\rho^-|$ we have the formula
$$|\widehat{K}_\rho^-|=|\widehat{K}_\rho|-|\widehat{K}|=C_{K_\rho}-C_K.$$
Thus, as a consequence of  Proposition \ref{prop-restrict} we see that
\begin{equation}\label{eq-O}
O_1+2O_2=|\widehat{K}_\rho^-|=C_{K_{\rho}}-C_K\quad\text{and}\quad 2O_1+O_2=|\widehat{G}_\rho^-|=C_{G_\rho}-C_G,
\end{equation}
As in the $\Pin^c$-case, these equations have unique solutions and, using 
 Corollary  \ref{cor-K}, we obtain
 
\begin{proposition}\label{prop-finite}
Let $\rho:G\to \O(V)$ be a linear action of the finite group $G$ on $V$. Then 
 \begin{equation}\label{eq-oriented}
\rank(K_G^*(V))=\left\{\begin{matrix} C_{G_\rho}-C_G& \text{if $*+\dim(V)=0 \mod 2$}\\
0 &\text{if $*+\dim(V)=1 \mod 2$}\end{matrix}\right\}
\end{equation}
if the action is orientation preserving, and 
 \begin{equation}\label{eq-nonoriented}
\rank(K_G^*(V))=\left\{\begin{matrix} \frac{1}{3}\big(2(C_{K_\rho}-C_K)-(C_{G_\rho}-C_G)\big)& \text{if $*+\dim(V)=0 \mod 2$}\\
\frac{1}{3}\big(2(C_{G_\rho}-C_G)-(C_{K_\rho}-C_K)\big)
  &\text{if $*+\dim(V)=1 \mod 2$}\end{matrix}\right\}
\end{equation}
if the action is not orientation preserving.
\end{proposition}

\begin{remark}\label{rem-Karoubi} Karoubi shows in \cite{Kar} that 
for any linear action $\rho:G\to \O(V)$ of a finite group $G$ 
the ranks of $K_G^0(V)$ and $K_G^1(V)$ can alternatively
be computed as follows: For any conjugacy class $C^g$ in $G$ let 
$V^g$ denote the fixed-point set of $\rho(g)$  in $V$. This space is 
$\rho(h)$-invariant for any $h$ in the centralizer $C_g$ of $g$,
 and therefore $C_g$ acts 
linearly on $V^g$ for all $g$ in $G$.
With these facts in mind, Karoubi denotes a conjugacy class $C^g$ \emph{oriented}, if the
action of $C_g$ on $V^g$ is oriented and $C_g$ is called \emph{even} (resp. \emph{odd})
if $\dim(V^g)$ is even (resp.  odd). He then shows in  \cite[Theorem 1.8]{Kar} that 
the rank of $K_G^0(V)$ (resp. $K_G^1(V)$) equals the number of oriented 
even (resp. odd) conjugacy classes in $G$. 
It seems not  obvious to us that Karoubi's result coincides with 
the result given in Proposition \ref{prop-finite} above.
\end{remark}

In what follows we want to apply our results  to give an alternative proof
of the formulas for the ranks of $K_{S_n}^0(\RR^n)$ and $K_{S_n}^1(\RR^n)$, as given by Karoubi in 
\cite[Corollary 1.9]{Kar}, where the symmetric group $S_n$, for $n\geq 2$,
acts on $\RR^n$ by permuting the standard orthonormal base $\{e_1,\dots, e_n\}$.
Let $\rho:S_n\to \O(n)$ denote the corresponding homomorphism.
 It is clear that the inverse image of $\SO(n)$ is the 
alternating group $A_n$. To simplify notation, we shall write $\tilde{S}_n$ and 
$\tilde{A}_n$ for the groups $(S_n)_\rho$ and $(A_n)_{\rho}$.
Thus it follows from Proposition \ref{prop-finite} that 
 \begin{equation}\label{eq-Sn}
\rank(K_{S_n}^*(\RR^n))=\left\{\begin{matrix}\frac{1}{3}\big(2(C_{\tilde{A}_n}-C_{A_n})-(C_{\tilde{S}_n}-C_{S_n})\big)  & \text{if $*+n=0 \mod 2$}\\
\frac{1}{3}\big(2(C_{\tilde{S}_n}-C_{S_n})-
(C_{\tilde{A}_n}-C_{A_n})\big) &\text{if $*+n=1 \mod 2$}\end{matrix}\right\}
\end{equation}
Thus, to get explicit formulas we need to compute the 
numbers $C_{\tilde{S}_n}-C_{S_n}$ and $C_{\tilde{A}_n}-C_{A_n}$. 
For this we use the following general observations:
If 
$$
\begin{CD} 1 @>>> \ZZ_2 @>>>\tilde{G} @>q>> G@>>>1
\end{CD}
$$
is any central extension of $G$ by $\ZZ_2$, the inverse image $q^{-1}(C^g)$ of a conjugacy class 
$C^g$ in $G$ is either  a conjugacy class in $\tilde{G}$ itself, or it decomposes into two disjoint 
 conjugacy classes of the same length in $\tilde{G}$. If $t\in \tilde{G}$ such that $q(t)=g$, the
  second possibility happens if and only if $t$ is not conjugate to $-t$ in $\tilde{G}$ (e.g., see \cite[Theorem 3.6]{HH}).
 Thus, if $C_G^{\dec}$ denotes the number of conjugacy classes in $G$ which decompose in $\tilde{G}$, the number $C_{\tilde{G}}$ of conjugacy classes in $\tilde{G}$ is equal to
  $C_G+C_G^{\dec}$, and hence
  $$C_{\tilde{G}}-C_G=C_G^{\dec}.$$
Now, for the groups $G=S_n$ and $K=A_n$ the numbers $C_{S_n}^{\dec}$ and $C_{A_n}^{\dec}$
have been computed explicitely  in \cite[Theorem 3.8 and Corollary 3.10]{HH}:

\begin{proposition}\label{prop-HH}
For each $n\geq 2$ let $a_n$ (resp. $b_n$) 
denote the number of all finite tupels of natural numbers
$(\lambda_1,\dots,\lambda_m)$ such that $1\leq \lambda_1<\lambda_2<\cdots <\lambda_m$,
$\sum_{i=1}^m\lambda_i=n$, and  the number of even entries $\lambda_i$ is even (resp. odd).
Then
$$C_{{S}_n}^{\dec}=a_n+2b_n\quad\text{and}\quad C_{{A}_n}^{\dec}=2a_n+b_n.$$
\end{proposition}

We should note that the definitions of $\tilde{S}_n$ and $\tilde{A}_n$ considered in \cite{HH}
are slightly different from ours, but a study of the proofs of Theorem 3.8 and Corollary 3.10 in that
paper shows that the arguments apply step by step to our situation. As a consequence we get

\begin{corollary}\label{cor-Sn}
Let $\rho:S_n\to\O(n)$ be as above. Then 
$$K_{S_n}^*(\RR^n)=\left\{\begin{matrix} \ZZ^{a_n} &\text{if $*+n=0 \mod 2$}\\
\ZZ^{b_n}&\text{if $*+n=1 \mod 2$}\end{matrix}\right\}.$$
\end{corollary}

\begin{remark}\label{rem-numbers} In \cite[Corollary 1.9]{Kar}, Karoubi gives  the formulas
$$K_{S_n}^0(\RR^n)=\ZZ^{p_n}\quad\text{and}\quad K_{S_n}^1(\RR^n)=\ZZ^{i_n}$$
where $p_n$ (resp. $i_n$) denotes the number of partitions 
$n = \sum_{i=1}^{m} \lambda_i$ 
with $1\leq \lambda_1<\cdots <\lambda_m$ and 
$m=2k$ even (resp. $m=2k+1$ odd).
One checks that the numbers $p_n$ and $i_n$ are related to the numbers $a_n$ and $b_n$,
as defined in the corollary above, by the equations:
$$\begin{array}{cc}
a_{2n+1}=i_{2n+1} & b_{2n+1}=p_{2n+1}\\
a_{2n}=p_{2n}& b_{2n}=i_{2n},
\end{array}
$$
and hence Karoubi's formula coincides with ours.
We give the argument for the equation
$a_{2n+1}=i_{2n+1}$, the other equations can be shown similarly.
So let $n\in \NN$ and let $(\lambda_1,\ldots, \lambda_m)$ be a partition of $2n+1$ as
in the definition of $a_{2n+1}$, i.e., there is an even number $2r$ of even entries $\lambda_i$ 
in this partition.  Since $\lambda_1+\cdots+\lambda_m=2n+1$ is odd, it follows that there 
is an odd number of odd entries $\lambda_j$ in the partition. Thus $m=2k+1$ is odd.
Conversely, if $m=2k+1$ is odd, the fact that  $\lambda_1+\cdots+\lambda_m=2n+1$ is odd
implies that the number $l$ of odd entries $\lambda_j$ must be odd, and then the number 
$m-l$ of even entries must be even.

If we restrict the action of $S_n$ to the alternating group $A_n$, we obtain the formulas
$$K_{A_n}^*(\RR^n)=\left\{\begin{matrix} \ZZ^{2a_n+b_n}&\text{if $*+n= 0\mod 2$}\\
\{0\}& \text{if $*+n=1\mod 2$}\end{matrix}\right\},$$
since this action is orientation preserving.
\end{remark}

\section{Actions of $\O(n)$}\label{sec-On}

In this section we want to study the $K$-theory groups $K_{\O(n)}^*(V)$ for linear actions 
$\rho:\O(n)\to \O(V)$ of the orthogonal group $\O(n)$ on an arbitrary real vector space $V$.
We are in particular interested in the canonical action of $\O(n)$ on $V=\RR^n$ and in the
action of $\O(n)$ on the space $V_n$ of all symmetric matrices in $M_n(\RR)$, with action given by
conjugation. The study of the latter will allow to compute explicitly the $K$-theory groups of 
the reduced group C*-algebra $C_r^*(\GL(n,\RR))$ of the general linear group $\GL(n,\RR)$ via the 
positive solution of the Connes-Kasparov conjecture.

 Recall that $\O(n)\cong \SO(n)\times \ZZ_2$ if $n$ is odd,
with $-I\in \O(n)$ the generator of $\ZZ_2$ (in what follows we denote by $I$ the unit matrix in $\O(n)\subseteq M_n(\RR)$ and we denote by $1$ the unit in $\Spin(n)\subseteq \Cl_\RR(n)$). If $n$ is even, we have
 $\O(n)\cong \SO(n)\rtimes \ZZ_2$, 
the semi-direct product of $\SO(n)$ with $\ZZ_2$, 
where the generator of $\ZZ_2$ can be chosen to be the 
matrix $g:=\diag(-1,1,\ldots,1)\in \O(n)$ acting on $\SO(n)$ by conjugation.

Given a representation $\rho:\O(n)\to \O(V)$, we need to describe the group $\O(n)_\rho$ and 
its representations. For this we start by describing all possible central extensions
of $\O(n)$ by $\ZZ_2$. Indeed, we shall see below that for any fixed $n\geq 2$ there are precisely four such extensions
$$\begin{CD}
1 @>>> \ZZ_2 @>>> G_i^n @>q>> \O(n) @>>> 1,
\end{CD}
$$
$i=0,\ldots, 3$. To describe them, 
we let  $K_i^n$ denotes the inverse image of $\SO(n)$ in $G_i^n$
for $i=0,\ldots, 3$. This is a central extension of $\SO(n)$ by $\ZZ_2$ and therefore 
the $K_i^n$ are either isomorphic to the trivial extension $\SO(n)\times \ZZ_2$ or the 
nontrivial extension $\Spin(n)$. Using this, the extensions $G_i^n$, $2\leq n\in \NN$, $i=0,\ldots,3$ are given as follows:\\
{\bf If $n=2m+1$ is odd}, then
\begin{enumerate}
\item[(O1)] there are two extensions $G_0^n$ and $G_1^n$ such that $K_0^n=K_1^n=\SO(n)\times \ZZ_2$:
the trivial extension $G_0^n=\O(n)\times \ZZ_2$ and the non-trivial extension $G_1^n=\SO(n)\times \ZZ_4$, with central subgroup $\ZZ_2$ being the order-two subgroup of $\ZZ_4$.
\item[(O2)] There are two extensions $G_2^n, G_3^n$ such that $K_2^n=K_3^n=\Spin(n)$.
To characterize them let $x\in G_i^n$ such that  $q(x)=-I\in \O(n)$. Then 
$x^2=1$ for $x\in G_2^n$ and $x^2=-1$  for $x\in G_3^n$. We then have $G_2^n\cong \Spin(n)\times \ZZ_2$ with $x$ a generator for $\ZZ_2$ and $G_3^n\cong (\Spin(n)\times \ZZ_4)/\ZZ_2$ with respect to the diagonal embedding of $\ZZ_2$ into $\Spin(n)\times \ZZ_4$. The central subgroup $\ZZ_2$ is 
given by (the image of)  the order-two subgroup $\{\pm 1\}\subseteq \Spin(n)$.
\end{enumerate}
{\bf If $n=2m$ is even}, then
\begin{enumerate}
\item[(E1)] there are two extensions $G_0^n$ and $G_1^n$ such that $K_0^n=K_1^n=\SO(n)\times \ZZ_2$:
the trivial extension $G_0^n=\O(n)\times \ZZ_2$ and the non-trivial extension 
$G_1^n=\SO(n)\rtimes \ZZ_4$, with action of $\ZZ_4$ on $\SO(n)$ given on the generator 
by conjugation with $g=\diag(-1,1,\ldots,1)$,
and the central subgroup $\ZZ_2$ of $\SO(n)\rtimes \ZZ_4$ 
is given by the order-two subgroup of $\ZZ_4$.
\item[(E2)] There are two extensions $G_2^n, G_3^n$ such that $K_2^n=K_3^n=\Spin(n)$.
If $x\in G_i^n$ such that  $q(x)=\diag(-1,1,\ldots,1)$, then $x^2=1$ in case $x\in G_2^n$ and $x^2=-1$
in case $x\in G_3^n$. Then $G_2^n\cong \Spin(n)\rtimes \ZZ_2$ with $\ZZ_2$ generated by $x$
and $G_3^n\cong (\Spin(n)\rtimes\ZZ_4)/\ZZ_2$ where $\ZZ_4$ is generated by an element $\tilde{x}
\in \Spin(n)\rtimes \ZZ_4$ which acts on $\Spin(n)$ by conjugation with $x$, and $\ZZ_2$ is embedded diagonally into $\Spin(n)\rtimes\ZZ_4$ as in the odd case. The central copy of $\ZZ_2$ is given by (the image of) 
the order-two subgroup $\{\pm 1\}\subseteq \Spin(n)$.
\end{enumerate}

Although we are convinced that this description of the central extensions of $\O(n)$ by $\ZZ_2$ 
is well-known, we give a proof since we didn't find a direct reference:

\begin{proposition}\label{prop-extOn}
For any fixed $n\geq 2$ the above described extensions are, up to isomorphism of extensions, the 
only central extensions of $\O(n)$ by $\ZZ_2$.
\end{proposition}
\begin{proof} Recall first that the set of isomorphism classes of central extensions of any given group
$H$
by $\ZZ_2$ forms a group $\mathcal E(H,\ZZ_2)$ with group operation given as follows: if
$$1\to \ZZ_2\to G\stackrel{q}{\to} H\to 1\quad\text{and}\quad
1\to \ZZ_2\to G'\stackrel{q'}{\to} H\to 1$$ are  central extensions of $H$ by $\ZZ_2$,
then the product is given by the (isomorphism class) of the extension
\begin{equation}\label{ext-product}
1 \longrightarrow\ZZ_2 \stackrel{\iota}{\longrightarrow}  G*G'\stackrel{q''}{\longrightarrow} H\longrightarrow 1,
\end{equation}
where $G*G'=\{(x,x')\in G\times G': q(x)=q'(x')\}/\ZZ_2$ 
with respect to the diagonal embedding of $\ZZ_2$  into 
$G\times G'$. The central copy of $\ZZ_2$ in $G*G'$ can be taken as the image in $G*G'$
of the central copy of $\ZZ_2$ in either $G$ or $G'$.
%
It is well known that $\E(\SO(n),\ZZ_2)\cong\ZZ_2$ with non-trivial element given by $\Spin(n)$. 
%
%

Suppose now that $1\to \ZZ_2\to G\stackrel{q}{\to} \O(n)\to 1$ represents an element in 
$\E(\O(n),\ZZ_2)$. It restricts to a representative $1\to \ZZ_2\to K\to \SO(n)\to 1$ in $\E(\SO(n),\ZZ_2)$
with $K:=q^{-1}(\SO(n))$. This restriction procedure
induces a homomorphism of 
$\E(\O(n),\ZZ_2)$ to $\E(\SO(n),\ZZ_2)$. 
Therefore, given any fixed extension $G$ which restricts to $K=\Spin(n)$, then all other extension which restrict to $\Spin(n)$ are given as products $G*G'$ where $G'$ is an extension 
which restricts to $K'=\SO(n)\times \ZZ_2$. In particular, if we can show that there are only two extensions which restrict to $\SO(n)\times \ZZ_2$, then there are also only two extensions which restrict 
to $\Spin(n)$. Since the ones given in the above list are obviously non-isomorphic (as extensions),
the list must be complete.

So let $G'$ be any extension 
which restricts to $K'=\SO(n)\times \ZZ^2$. We show that it equals $G_0^n$ or  $G_1^n$ described above.
Suppose first that $n=2m$ is even.
Choose $x\in G'$  such that $q(x)=g:=\diag(-1,1,\ldots,1)\in \O(n)$.
Then $q(x^2)=I$  and hence $x^2=\pm 1$. We claim that 
$$x(h, \epsilon)x^{-1}=(ghg^{-1},\epsilon)$$
for all $(h,\epsilon)\in \SO(n)\times \ZZ_2$. 
To see this, note first that 
$q\big(x(h, \epsilon)x^{-1}\big)=g\big(q(h, \epsilon)\big)g^{-1}=ghg^{-1}$, which implies that
$x(h, \epsilon)x^{-1}=(ghg^{-1},\epsilon')$ for some $\epsilon'\in \ZZ_2$. To see
that $\epsilon'=\epsilon$ we simply observe that the map $\SO(n)\to \ZZ_2$ which sends
$h\in \SO(n)$  to the projection 
of $x(h,1)x^{-1}$ to $\ZZ_2$ is a continuous group homomorphism, and hence trivial.
This implies $1'=1$ and then  also $(-1)'=-1$.

If $x^2=1$, it follows that $G= (\SO(n)\times\ZZ_2)\rtimes \lk x\rk=(\SO(n)\rtimes\lk g\rk)\times \ZZ_2
=\O(n)\times \ZZ_2=G_0^n$.
If $x^2=-1$, we obtain a surjective homomorphism 
$\varphi:(\SO(n)\times\ZZ_2)\rtimes  \lk x\rk\to G$ given by $\varphi\big(x^j, (g, \pm 1)\big)=x^j(g,\pm 1)$
with kernel generated by the order-two element $(x^2,-1)$. This implies $G=G_1^n$.

A similar but easier argument applies in the case where $n=2m+1$ is odd. We omit the details for this.
\end{proof}

\begin{remark} As we saw in  Section \ref{sec-prel} (e.g. see (\ref{ex-even}) and (\ref{ex-odd}))
the group  $\Pin(n)$ is a central extension of $\O(n)$ by $\ZZ_2$ 
which restricts to $\Spin(n)$.  Thus the other such extension is given by the product 
$\Pin(n)*G_1^n$. To see whether $\Pin(n)$ is the group $G_2^n$ or the group $G_3^n$,
we need to identify an inverse image $x\in \Pin(n)$ of the matrix $-I$ if $n=2m+1$ is odd, or of
 $g=\diag(-1,1,\dots,1)\in \O(n)$ if $n=2m$ is even.

Indeed, if $e_1,\ldots, e_n$ denotes the standard orthonormal base of $\RR^n$, then it follows from the basic relations in in $\Cl_\RR(n)$ that
$x=\pm e_1\cdots e_n$ if $n=2m+1$ is odd (with respect to extension (\ref{ex-odd})) and 
$x= \pm e_2\cdots e_n$ if $n=2m$ is even.
In the first case we get 
$$x^2=(-1)^{\frac{n(n+1)}{2}}=(-1)^{(m+1)(2m+1)}=\left\{\begin{matrix} -1&\text{if $m$ is even}\\
1&\text{if $m$ is odd}\end{matrix}\right\},$$
and in the second case we  get
 $$x^2=(-1)^{\frac{n(n-1)}{2}}=(-1)^{m(2m-1)}=\left\{\begin{matrix} 1&\text{if $m$ is even}\\
-1&\text{if $m$ is odd}\end{matrix}\right\}.$$
\end{remark}

In what follows we need to understand the conjugation action of $G_i^n$
on $\widehat{K_i^n}$. Note that in all cases we can identify $G_i^n/K_i^n$ 
with $\ZZ_2$. 

\begin{lemma}\label{lem-conjugation}
Let $n\geq 2$,  let $1\to \ZZ_2\to G\stackrel{q}{\to} \O(n)\to 1$ be any central extension 
of $\O(n)$ by $\ZZ_2$ and let $K=q^{-1}(\SO(n))$. Then 
%
\begin{enumerate}
\item If $n$ is odd, the conjugation action of $G/K$  on $\widehat{K}$ is trivial.
\item If $n$ is even, and if $K=\SO(n)\times\ZZ_2$,
 the action of  $G/K$ on $\widehat{K}=\widehat{\SO(n)}\times\widehat{\ZZ}_2$ is given by 
the conjugation action of $\O(n)/\SO(n)$ on the first factor and the trivial action on the second.
If  $K=\Spin(n)$, the action of $G/K$ on $\widehat{\Spin(n)}$ coincides with the conjugation action of $\Pin(n)/\Spin(n)$ on $\widehat{\Spin(n)}$.
\end{enumerate}
\end{lemma}
\begin{proof} 
The first assertion follows directly from the description of the groups $G_i^n$ in case where $n$ is odd.
So assume now that $n$ is even and $K=\SO(n)\times \ZZ_2$. Let $x\in G$ with $q(x)=g:=\diag(-1,1\ldots,1)$. It is shown in the proof of Proposition \ref{prop-extOn} that the conjugation 
action of $x$ on $\SO(n)\times \ZZ_2$ is given by conjugation with $g\in \O(n)$ in the first factor and 
the trivial action in the second factor. This proves the first assertion in (ii). 

So assume now that 
$K=\Spin(n)$. Then $G=\Pin(n)$ or $G=\Pin(n)*G_1^n$. The result is clear in the first case.
So let $G=\Pin(n)*G_1^n$. 
If $y=e_2\cdots e_n\in \Pin(n)$ and $x_1\in G_1^n$ with $q_1(x_1)=g$, then
$x=[y,x_1]$ is an inverse image of $g$ in $G$. The group $\Spin(n)$ then identifies with 
$K\subseteq G$  via the embedding
$$\varphi:\Spin(n)\to \Pin(n)*G_1^n; z\mapsto [z,(\Ad(z), 1)].$$
 Conjugating  $ [z,(\Ad(z), 1)]$ by $[y,x_1]$  provides  
 $[yzy^{-1}, x_1(\Ad(z),1)x_1^{-1}]=[yzy^{-1}, (g\Ad(z) g^{-1},1)]=\varphi(yzy^{-1})$,
 which  finishes the proof.
\end{proof}

Suppose now that $\rho:\O(n)\to \O(V)$ is any linear action of $\O(n)$ on a finite dimensional 
real vector space $V$. 
In all cases the group $\O(n)_\rho$ must be one of the groups $G_0^n,\ldots, G_3^n$.
If we fix $n\geq 2$, we have the following possibilities for the computation of the $K$-theory groups 
$K_{\O(n)}^*(V)$:\\
\\
{\bf The orientation preserving case:} 
If the action of $\O(n)$ on $V$ is orientation preserving, then Corollary \ref{cor-oriented}
implies that 
\begin{equation}\label{eq-oriented1}
K_{\O(n)}^*(V)\cong \left\{\begin{matrix}\oplus_{\sigma\in \widehat{G_i^n}^-}\ZZ&\text{if $*+\dim(V)$ is even,}\\
\{0\}& \text{if $*+\dim(V)$ is odd}\end{matrix}\right\}.
\end{equation}
Note that although we have 
four different possibilities $G_0^n,\ldots, G_3^n$ for the groups $\O(n)_\rho$, 
it follows from  Lemma \ref{lem-reps} that the cardinality of  $\widehat{G_i^n}^-$ is always countably infinite. Thus the isomorphism class  
of $K_{\O(n)}^*(V)$ only depends on the fact whether $\dim(V)$ is even or odd.
\\
\\
{\bf The non-orientation preserving case:} If the action of $\O(n)$ is not orientation preserving, the question which group  
out of the list $G_0^n,\ldots, G_3^n$ we get for $\O(n)_\rho$ becomes more interesting (at least if $n$ is even).
In fact, if $i\in \{0,\ldots, 3\}$ such that 
$\O(n)_\rho\cong G_i^n$ as central extension of $\O(n)$ by $\ZZ_2$, it follows from Corollary \ref{cor-K} 
that
\begin{equation}\label{eq-orbit}
K_{\O(n)}^*(V)\cong \left\{\begin{matrix}\oplus_{\sigma\in \mathcal O_1} \ZZ&\text{if $*+\dim(V)$ is odd,}\\
\oplus_{\{\sigma,g\sigma\}\in \mathcal O_2} \ZZ&\text{if $*+\dim(V)$ is even}\end{matrix}\right\},
\end{equation}
where $\mathcal O_1$ and $\mathcal O_2$ denote the numbers of orbits  
 of length one or  two in the set $\widehat{K}^-$ of negative representations of $K:=K_i^n$ under the
 conjugation action of $G_i^n$.
 So in order to get the general picture, we need to study the cardinalities of the sets 
 $\mathcal O_1$ and $\mathcal O_2$ in the four possible cases. We actually get different answers 
 depending on whether $n$ is even or odd:
\\
\\
{\bf The odd case $n=2m+1$:}  In this case 
 Lemma \ref{lem-conjugation} implies that the action of $G_i^n$ on $\widehat{K}^-$ is trivial 
 in all cases. Thus, from the above formula we get
 \begin{equation}\label{eq-odd}
 K_G^*(V)\cong \left\{\begin{matrix}\oplus_{\sigma\in \widehat{K}^-} \ZZ&\text{if $*+\dim(V)$ is odd}\\
\{0\}&\text{if $*+\dim(V)$ is even}\end{matrix}\right\}.
\end{equation}
As in the orientation preserving case, it follows from Lemma \ref{lem-reps} that
 the cardinality of $\widehat{K}^-$ is always countably infinite.
 \\
 \\
{\bf The even case $n=2m$:}  Let $G=\O(n)_\rho$ and $K=q^{-1}(\SO(n))$. 
If $K=\SO(n)\times\ZZ_2$
we get $\widehat{K}^-=\widehat{\SO(n)}\times\{\mu\}\cong \widehat{\SO(n)}$, where 
 $\mu$ denotes the non-trivial character of $\ZZ_2$, and it follows from Lemma \ref{lem-conjugation}
 that the action of $G/K$ on $\widehat{K}^-$  is given by the conjugation action  of $\O(n)/\SO(n)$ on $ \widehat{\SO(n)}\cong \widehat{K}$.
  Thus, the orbit sets $\mathcal O_1$ and $\mathcal O_2$ 
can be identified with the $\O(n)$-orbits of length one and two in $\widehat{\SO(n)}$ and the 
$K$-theory groups in the cases $\O(n)_\rho=G_0^n$ and $\O(n)_\rho=G_1^n$ are the same 
(up to isomorphism).\medskip
\\
In case $K=\Spin(n)$ it follows  from Lemma \ref{lem-conjugation} that  the
 action of $G/K$ on $\widehat{\Spin(n)}$ coincides with the
  conjugation action of $\Pin(n)/\Spin(n)$ on $\widehat{\Spin(n)}$.
Thus, to compute the sets $\mathcal O_1$ and $\mathcal O_2$ in the $K$-theory formula
(\ref{eq-orbit}), we may assume 
without loss of generality that $\O(n)_\rho=\Pin(n)$.

\medskip
In view of the above discussions, it is desirable to find an easy criterion 
for the group $K\subseteq \O(n)_\rho$ being isomorphic to $\SO(n)\times \ZZ_2$ or not.
It is clear that  this is the case if and only if 
the restriction $\rho:\SO(n)\to \SO(V)$
of  the given action $\rho:\O(n)\to\O(V)$ is {\em spinor} in the sense that  there exists a homomorphism
$\tilde\rho:\SO(n)\to \Spin(V)$ such that $\rho=q\circ \tilde{\rho}$ with $q:\Spin(V)\to \SO(V)$ the quotient map.

\begin{proposition}\label{prop-spin}
Let $n\geq 2$ and let $\rho:\SO(n)\to \SO(V)$ be any linear action of $\SO(n)$ on the finite dimensional real vector space $V$. Let $h:=\diag(-1,-1,1,\ldots, 1)\in \SO(n)$ and let $V^-=\{v\in V: \rho(h)v=-v\}$
denote the  eigenspace for the eigenvalue $-1$ of $\rho(h)$. Then $\rho$ is spinor if and only 
if $\dim(V^-)=4k$ for some $k\in \NN_0$.
\end{proposition}
\begin{proof} We first need to know that, given a representation $\rho:\SO(n)\to \SO(V)$, 
there always exists a representation $\hat{\rho}:\Spin(n)\to \Spin(V)$ such that the diagram
\begin{equation}\label{eq-rhodach}
\begin{CD}
\Spin(n) @>\hat\rho >>\Spin(V)\\
@V q=\Ad VV   @VVq=\Ad V\\
\SO(n)  @>>\rho> \SO(V)
\end{CD}
\end{equation}
commutes. In case $n>2$ this follows from the universal properties of the universal covering
$\Spin(n)$ of $\SO(n)$. In case $n=2$, the groups $\Spin(2)$ and $\SO(2)$ are both isomorphic to the 
circle group $\TT$ with covering map $q:\TT\to\TT; z\mapsto z^2$. The image $\rho(\SO(2))$
lies in a maximal torus $T\subseteq \SO(V)$ and there is a maximal Torus $\tilde{T}$ in 
$\Spin(V)$ which projects onto $T$ via a double covering map $q:\tilde{T}\to T$. Thus the 
problem reduces to the problem whether there exists a map $\hat{\rho}:\TT\to \tilde{T}$
which makes the diagram
$$
\begin{CD}
\TT @>\hat\rho >>\tilde{T}\\
@V qVV   @VVqV\\
\TT  @>>\rho> T
\end{CD}
$$
commute. It is straightforward to check that this is always possible.

It follows from (\ref{eq-rhodach}) that there exists a lift $\tilde\rho:\SO(n)\to \Spin(n)$ for $\rho$
if and only 
$$\{\pm 1\}=\ker\big(q:\Spin(n)\to \SO(n)\big)\subseteq \ker\hat\rho.$$
So we simply have to check whether $\hat\rho(-1)=1$ or not.
For this let $e_1,\ldots, e_n$ denote the standard orthonormal base of $\RR^n$.
Then the product $e_1e_2\in \Spin(n)$ projects onto $h$ and $(e_1e_2)^2=-1$.
Since $\rho(h)^2=\rho(h^2)=1$ we see that
$V$ decomposes into the orthogonal direct sum $V^+\oplus V^-$ 
with $V^+$ and $V^-$ the eigenspaces for $\pm1$ of $\rho(h)$. 
Since $\det(\rho(h))=1$, it follows that $l=\dim(V^-)$ is even.
If $V^-=\{0\}$ we have $\rho(h)=1$, which implies $\hat\rho(e_1e_2)=\pm1$ and hence
$\hat\rho(-1)=\hat\rho((e_1e_2)^2)=1$.

If $V^-\neq\{0\}$ let  $\{v_1,\ldots v_l\}$ be any orthonormal base for $V^-$. It then
follows from the relations
in $\Cl_\RR(n)$ that
the element
$$y:=v_1\cdots v_l\in \Spin(V)$$
projects onto $\rho(h)\in \SO(V)$. This implies that $\hat{\rho}(e_1e_2)=\pm y$ and hence
that 
$$\hat\rho(-1)=\hat\rho((e_1e_2)^2)=y^2=(v_1\cdots v_l)^2=(-1)^{\frac{l(l+1)}{2}}.$$
Since $l$ is even, we get $l=2m$ for some $m\in \NN$ and then
$$\hat\rho(-1)=(-1)^{m(2m+1)}=\left\{\begin{matrix} 1&\text{if $m$ is even}\\
-1& \text{if $m$ is odd}\end{matrix}\right\}.$$
This finishes the proof.
\end{proof}

The only problem  which now remains for  the general
computation of $K_{\O(n)}^*(V)$ is the problem of computing explicitly
the orbit sets $\mathcal O_1$ and $\mathcal O_2$ which appear in 
formula (\ref{eq-orbit}) in the case where $n$ is even (as observed above, we always have 
$\mathcal O_1$ countably infinite and $\mathcal O_2=\emptyset$ if $3\leq n=2m+1$ is odd).
In order to give the complete picture, we now state the general result, although we postpone 
the proof for the case $n>2$ to  \S\ref{sec-orbits} below:

\begin{theorem}\label{prop-general-even}
Suppose that $\rho:\O(n)\to \O(V)$ is a non-orientation preserving action 
of $\O(n)$ on $V$ with $n=2m$ even. Then the following are true:
\begin{enumerate}
\item If the restriction $\rho:\SO(n)\to \SO(V)$ is spinor, then 
$\mathcal O_1$ consists of a single point if $n=2$ and $\mathcal O_1$ is countably infinite 
if $n>2$. The set $\mathcal O_2$  is always countably infinite.
\item If $\rho:\SO(n)\to \SO(V)$ is not spinor, then $\mathcal O_1=\emptyset$ and 
$\mathcal O_2$ is countably infinite.
\end{enumerate}
\end{theorem}

Combining this result with (\ref{eq-orbit})  immediately gives
\begin{corollary}\label{cor-general-even}
Suppose that $\rho:\O(n)\to \O(V)$ is a non-orientation preserving action 
of $\O(n)$ on $V$ with $n=2m$ even. Then
$$
K_{\O(n)}^0\cong\oplus_{n\in \NN}\ZZ
\quad\text{and}\quad
K_{\O(n)}^1(V)\cong\left\{\begin{matrix}\ZZ, &\text{if $n=2$,}\\
\oplus_{n\in\NN}\ZZ, &\text{if $2<n=2m$}\end{matrix}\right\}
$$
if the restriction $\rho:\SO(n)\to\SO(V)$ is spinor. Otherwise we get
$$
K_{\O(n)}^0(V)\cong\oplus_{n\in \NN}\ZZ
\quad\text{and}\quad
K_{\O(n)}^1=\{0\}.
$$
\end{corollary}

This corollary together with the discussions on the odd case implies

\begin{theorem}\label{thm-On} Let $\O(n)$ act on $\RR^n$ by matrix multiplication. Then
$$K_{\O(n)}^0(\RR^n)=\oplus_{k\in \NN}\ZZ\quad\text{and}\quad K_{\O(n)}^1(\RR^n)=\{0\}$$
 for
all $n\in \NN$ with $n\geq 2$.
\end{theorem}
\begin{proof}
Since the action is not orientation preserving and the restriction of $\id:\O(n)\to\O(n)$ to $\SO(n)$
is not spinor (which is an easy consequence of Proposition \ref{prop-spin}),  the result follows from formula (\ref{eq-odd}) in case where $n$ is odd, and from Corollary \ref{cor-general-even} if $n$ is even.
\end{proof}

The case $n=2$ of Theorem \ref{prop-general-even} is quite easy and has to 
be done separately, since the 
general methods used for $n>2$ in \S\ref{sec-orbits} below will not apply to this case. 
So we do the  case $n=2$ now:

\begin{proof}[Proof of Theorem \ref{prop-general-even} in case $n=2$]
As usual, let $K=q^{-1}(\SO(2))$ denote the inverse image of $\SO(2)$ in $\O(2)_\rho$.
If $\rho:\SO(2)\to \SO(V)$ is spinor, we have $K=\SO(2)\times \ZZ_2$. Otherwise we have
$K=\Spin(2)$.
\medskip
\\
{\bf The case $K=\SO(2)\times\ZZ_2$:} It follows from Lemma \ref{lem-conjugation} that
in this case the sets $\mathcal O_1$ and $\mathcal O_2$ can be identified with the sets of $\O(2)$-orbits 
in $\widehat{\SO(2)}$ of length
one and two, respectively.
Writing 
$$\SO(2)=\left\{g_\alpha=\left(\begin{smallmatrix}\cos(\alpha)&\sin(\alpha)\\-\sin(\alpha)&\cos(\alpha)
\end{smallmatrix}\right): \alpha\in [0,2\pi)\right\}$$
we  have $\widehat{\SO(2)}=\{\chi_k: k\in \ZZ\}$ with $\chi_k(g_\alpha)=e^{ik\alpha}$.
The action of $\O(2)$ on $\widehat{\SO(2)}$ is given by conjugation with $g=\diag(-1,1)$.
Since $gg_\alpha g^{-1}=g_{-\alpha}$ we get $g\cdot \chi_k=\chi_{-k}$, which implies that
$$\mathcal O_1=\{\chi_0\}\quad\text{and}\quad \mathcal O_2=\{\{\chi_k,\chi_{-k}\}: k\in \NN\}.$$
{\bf  The case $K=\Spin(2)$:} In this case the sets $\mathcal O_1$ and $\mathcal O_2$ can be identified with the sets of $\Pin(2)$-orbits 
in $\widehat{\Spin(2)}^-$ of length
one and two, respectively.
Recall that $\Spin(2)$ can be described as
$$\Spin(2)=\{x(\alpha):=\cos(\alpha) 1+\sin(\alpha)e_1e_2: \alpha\in [0,2\pi)\}\subseteq \Cl_\RR(2).$$
Then $\widehat{\Spin(2)}=\{\chi_k:k\in \ZZ\}$ with
$\chi_k:\Spin(2)\to \TT$, $ \chi_k(x(\alpha))=e^{i k\alpha }$.
It follows that $\widehat{\Spin(2)}^-=\{\chi_{2m+1}: m\in \ZZ\}$.
The action of $\Pin(2)$ on $\widehat{\Spin(2)}$ is given by conjugation with 
$x=e_2$. A short computation shows that 
$$e_2x(\alpha)e_2^*=-e_2x(\alpha)e_2=x(-\alpha)$$
which implies that $x\cdot \chi_k=\chi_{-k}$ for all $k\in \ZZ$. We therefore get
$x\cdot\chi_{2m+1}=\chi_{-2m-1}\neq \chi_{2m+1}$ for all $\chi_{2m+1}\in \widehat{\Spin(2)}^-$.
Thus
$$\mathcal O_1=\emptyset\quad\text{and}\quad \mathcal O_2=\{\{\chi_{2m+1},\chi_{-2m+1}\}: m\in \NN\}.$$
\end{proof}

We close this section with another interesting application of our main results.
Recall that for a locally compact group $H$, the reduced 
group $C^*$-algebra $C_r^*(H)$ is the closure of $\lambda(L^1(H))\subseteq \B(L^2(H))$, where
$$\lambda:L^1(H)\to \B(L^2(H)); \lambda(f)\xi=f*\xi$$
denotes the left regular representation of $H$. If $H$ is almost connected, it follows from the 
positive solution of the Connes-Kasparov conjecture, that there is a (more or less) canonical 
isomorphism $K^*(C_r^*(H))\cong K_G^*(V)$, where $G\subseteq H$ denotes the maximal 
compact subgroup of $H$
and  $V=T_{eG}(H/G)$ denotes the tangent space of $H/G$ at the 
trivial coset $eG=G$. The action of $G$ on $V$ is given by the 
differential of the left translation action  of $G$ on the manifold $H/G$ (see \cite[\S 7]{CEN}).

In case where $H=\GL(n,\RR)$, the
 maximal compact subgroup 
is $\O(n)$. If $V_n=\{A\in M(n,\RR): A=A^t\}$ denotes the space of symmetric matrices in $M_n(\RR)$, we have the well-known diffeomorphism
$$ V_n\times \O(n)\to \GL(n,\RR); \;(A,g)\mapsto \exp(A)g,$$
with $\exp(A)=\sum_{n=0}^\infty\frac{1}{n!}A^n$ the usual exponential map. 
Composing $\exp$ with the quotient map $\GL(n,\RR)\to \GL(n,\RR)/\O(n)$
provides a diffeomorphism $\widetilde{\exp}:V_n\to \GL(n,\RR)/\O(n)$.
We then get $\widetilde{\exp}(gAg^{-1})=g\cdot\widetilde{\exp}(A)$ and it follows 
from the above discussion that
\begin{equation}\label{eq-GLn}
K_*\big(C_r^*(\GL(n,\RR))\big)\cong K_{\O(n)}^*(V_n)
\end{equation}
for all $n\geq 2$, where the action of $\O(n)$ on $V_n$ is given by 
the representation 
$$\rho:\O(n)\to\O(V_n);\;\rho(g)A=gAg^{-1}$$
for all $g\in \O(n)$, $A\in V_n\subseteq M(n,\RR)$.

\begin{lemma}\label{lem-Vn}
Let $\rho:\O(n)\to \O(V_n)$ be as above. Then
\begin{enumerate}
\item $\rho$ is orientation preserving if and only if $n$ is odd.
\item The restriction $\rho: \SO(n)\to\SO(V_n)$ is spinor if and only if $n$ is even.
\end{enumerate}
\end{lemma}
\begin{proof} Let $\{E_{ij}: 1\leq i\leq j\leq n\}$ denote the standard 
basis of $V_n$, i.e., $E_{ij}$ has entry $1$ at the $ij$-th and $ji$-th place, and 
$0$ entries  everywhere else.
Conjugation with $g=\diag(-1,1,\ldots,1)\in \O(n)$ maps $E_{1j}$ to $-E_{1j}$ for all $j>1$ 
and fixes all other $E_{ij}$'s. It thus follows that $\det(\rho(g))=(-1)^{n-1}$, which shows that $\rho$ is 
orientation preserving if and only if $n$ is odd.

For the proof of (ii) we use  Proposition \ref{prop-spin}: let
$h=\diag(-1,-1,1\ldots,1)\in \SO(n)$.
Then conjugation with $h$ maps $E_{ij}$ to  $-E_{ij}$ for all $i=1,2$ and $j>2$ and fixes all other
$E_{ij}$. Thus $\{E_{ij}: i=1,2, j>2\}$ forms a base for $V_n^-$, the eigenspace of $\rho(h)$ 
for the eigenvalue $-1$. We therefore get $l:=\dim(V_n^-)=2(n-2)$.
This is a multiple of $4$ if and only if  $n$ is even.
\end{proof}

\begin{theorem}\label{thm-GLn}
If $n=2m+1$ is odd, then
$$ K_*\big(C_r^*(\GL(n,\RR))\big)\cong\left\{\begin{matrix} \oplus_{n\in \NN} \ZZ&\text{if $*+m$ is odd}\\
\{0\}&\text{if $*+m$ is even}\end{matrix}\right\}.$$
If $n=2m\geq 4$ is even, we get 
$$
 K_0(C_r^*(\GL(n,\RR))\cong \oplus_{n\in \NN} \ZZ\cong K_1\big(C_r^*(\GL(n,\RR))\big),$$
 and for $n=2$ we get
 $$ K_0\big(C_r^*(\GL(2,\RR))\big)\cong\ZZ\quad\text{and}\quad
 K_1\big(C_r^*(\GL(2,\RR))\big)\cong \oplus_{n\in\NN}\ZZ.$$
 \end{theorem}
 \begin{proof}
We use formula (\ref{eq-GLn}).  If $n=2m+1$ is odd, the result  then follows directly from formula
 (\ref{eq-oriented1}) together with 
 Lemma \ref{lem-Vn} above and the fact that $\dim(V_n)=\frac{n(n+1)}{2}=(2m+1)(m+1)$ is 
 even if and only if $m$ is odd. 
 
If $n$ is even, it follows from Lemma \ref{lem-Vn} above
 that the action of $\O(n)$ on $V_n$ 
is not orientation preserving and the restriction of $\rho$ to $\SO(n)$ is spinor. Thus the result 
 follows from Corollary \ref{cor-general-even}.
 \end{proof}

\section{Orbits in $\widehat{\Spin(n)}^-$ and $\widehat{\SO(n)}$}\label{sec-orbits}

In this section we want to provide the theoretical background to complete the proof of 
Theorem \ref{prop-general-even}. We need to compute the cardinalities for the orbit sets
$\mathcal O_1$ and $\mathcal O_2$  in $\widehat{\Spin(m)}^-$ under the conjugation action 
of $\Pin(m)$ and similarly for
the conjugation action of $\O(n)$ on $\widehat{\SO(n)}$. 

To solve this problem, we need
some background on the representation theory of  a connected compact Lie group $G$.
We use  \cite[Chapter VI]{BtD} as a general reference. 
Let $T$ denote a maximal torus in $G$ and let $\frak t$ denote its Lie algebra.
Let $I^*\subseteq \frak t^*$ denote the set of integral weights on $T$, i.e.,
the set of linear functionals $\lambda:\frak t\to \RR$ which vanish on the kernel
of $\exp:\frak t\to T$. There is a one to one correspondence between 
$I^*$ and $\widehat{T}$ given by sending an integral weight $\lambda$ to the character 
$e_\lambda:T\to \TT$ defined by $e_{\lambda}(\exp(t))=e^{2\pi i \lambda(t)}$ for all $t\in \frak t$.

Let $\bar{C}$ denote the closure of a fundamental Weyl chamber $C$ in $\frak t^*$
and let $\theta_1,\ldots, \theta_l\in I^*\cap \bar{C}$ be the corresponding set of positive roots.
In particuar,  $\theta_1,\ldots, \theta_l$ is a base of $\frak t^*$ and
$ \bar{C}=\{\sum_{i=1}^l a_i\theta_i: a_i\geq 0\}$.
There is a natural order on $\bar{C}$ given by $\lambda\leq\eta\Leftrightarrow \eta-\lambda\in \bar{C}$.
Let $W$ be the Weyl-group of  $G$, i.e., the group of automorphisms of $T$ induced from inner
automorphisms of $G$. Then $W$ acts canonically on $T$, $\frak t$, $\frak t^*$ and $I^*$. 

For any finite dimensional complex 
 representation $\tau$ of $G$ the equivalence class of $\tau$ is uniquely determined by its character 
 $\chi_{\tau}:=\tr\tau$, which is constant on conjugacy classes in $G$.
 A virtual character is a linear combination of such characters with integer coefficients. 
 The set $R(G)$ of all virtual characters of $G$ is called the representation ring of $G$.
 It is actually a subring of the ring of continuous functions on $G$.
Every element in $R(G)$ can be written as a (integer) linear combination of irreducible characters, i.e., the characters corresponding to 
 irreducible representations of $G$. 
 Since the restriction $\tau|_T$ of a representation $\tau$ of $G$ is invariant under conjugation with elements in $W$ (up to equivalence), the restriction of its character $\chi_\tau$ to $T$ is 
 conjugation invariant, and hence lies in the set $R(T)^W$ of {\em symmetric} (i.e., $W$-invariant)
 virtual characters of $T$.  By  \cite[Chapter VI, Proposition (2.1)]{BtD}
the restriction map 
 $$\res:R(G)\to R(T)^W; \chi\mapsto \chi|_T$$
 is an isomorphism of rings. 
 Now, for any $\lambda\in I^*$ we let $W\lambda=\{w\cdot\lambda: w\in W\}$ denote the $W$-orbit of $\lambda$ in $I^*$ and let
 $$S(\lambda)=\sum_{\xi\in W\lambda}e_{\xi}$$
 denote the {\em symmetrized} character corresponding to $\lambda$. 
 Combining  \cite[Chapter VI, Theorem (1.7) and Proposition (2.6)]{BtD} we get the following version
 of Weyl's character formula:

 \begin{theorem}\label{thm-Weyl}
 For each irreducible representation $\tau$ of $G$ there exists a unique decomposition
 $$\chi_{\tau}|_T=S(\lambda)+\sum_{i=1}^k l_i S(\lambda_i)$$
 with pairwise different $\lambda, \lambda_1,\ldots, \lambda_k\in I^*\cap\bar{K}$, $l_1,\ldots, l_k\in \ZZ$ and $\lambda_i< \lambda$
 for all $1\leq i\leq k$. We call $\lambda\in I^*\cap\bar{C}$ the {\em highest weight} of the 
 representation $\tau$. The map which assigns $\tau$ to its highest weight $\lambda$ 
 induces a bijection  between  $\widehat{G}$ and $I^*\cap \bar{C}$.
 \end{theorem}
 
 In what follows we shall denote by $\chi_\lambda\in R(G)$ the character of the 
 irreducible representation $\tau$ with highest weight $\lambda$. If $\gamma, \lambda$ are weights in 
 $I^*\cap\bar{C}$, then so is $\gamma+\lambda$ and there is a corresponding irreducible 
 character $\chi_{\gamma+\lambda}$ of $G$. By \cite[Chapter VI, (2.8)]{BtD} we have
 
 \begin{lemma}\label{lem-product}
 For all $\gamma, \lambda\in I^*\cap\bar{C}$ there is a unique decomposition
 $$\chi_\gamma\cdot\chi_{\lambda}=\chi_{\gamma+\lambda}+\sum_{\mu}l_\mu\chi_{\mu},$$
 where $\mu$ runs through $\{\mu\in I^*\cap\bar{C}: \mu<\lambda+\gamma\}$ and 
 $0\leq l_\mu\in \ZZ$.
 \end{lemma}
 
 A set $\{\lambda_1,\ldots, \lambda_k\}$ of integral weights in $I^*\cap \bar{C}$ is called 
 a {\em fundamental system}, if the map
 $$\varphi: \NN_0^k\to I^*\cap\bar{C}; \;\varphi(l_1,\ldots, l_k)=\sum_{i=1}^k l_i\lambda_i$$
 is an ordered bijection with respect to the standard order on $\NN_0^k$. 
 The corresponding irreducible representations $\tau_1,\ldots,\tau_k$ are then called 
 a system of fundamental representations of $G$.
 By \cite[Chapter VI (2.10) and (2.11)]{BtD} we have  
 
 \begin{theorem}\label{thm-fund}
 Suppose that $G$ is a connected and simply connected compact Lie group. Then 
 there exists a fundamental system $\{\lambda_1,\ldots, \lambda_k\}$ in $I^*\cap \bar{C}$ and 
 there is a ring isomorphism 
 $\psi: \ZZ[X_1,\ldots, X_k]\to R(G)$ which sends $X_i$ to $\chi_{\lambda_i}$.
 \end{theorem}
 
  Combining these results, we get

\begin{lemma}\label{lem-polydeco} 
Let $G$ be a connected and simply connected Lie group and let $\{\lambda_1,\ldots, \lambda_k\}$
be a fundamental system in $I^*\cap\bar{C}$. Let $\gamma=\sum_{i=1}^k l_i\lambda_i$ 
be any given weight in $I^*\cap \bar{C}$ and let $\chi_\gamma$ be the corresponding irreducible character of $G$. Then there exists a unique polynomial $P\in \ZZ[X_1,\ldots, X_k]$ of order less than 
$l:=l_1+l_2+\ldots+l_k$ such that 
$$\chi_{\gamma}=\prod_{i=1}^k\chi_{\lambda_i}^{l_i}+P(\chi_{\lambda_1},\ldots, \chi_{\lambda_k}).$$
\end{lemma} 
\begin{proof}
Uniqueness is a direct consequence of Theorem \ref{thm-fund} above.
For existence, we give a proof by induction on the sum $l=l_1+l_2+\ldots+l_k$ corresponding to $\gamma$, which 
we call the {\em order} of $\gamma$.
If $l=0$, then  $\chi_{\gamma}\equiv 1$ is the character of the trivial representation and 
the formula is true with $P=0$ (we use the convention that the order of the zero-polynomial is $-\infty$).
Suppose now that for given $l> 0$ the lemma is true for all $m<l$. Let $\gamma\in I^*\cap\bar{C}$ 
with order $l$, $\gamma=\sum_{i=1}^kl_i\chi_{\lambda_i}$. Without loss of generality we may assume
that $l_1>0$. By Lemma \ref{lem-product} we have
$$\chi_{\gamma}=\chi_{\lambda_1}\chi_{\gamma-\lambda_1}-\sum_{\mu<\gamma}l_{\mu}\chi_{\mu}$$
 for suitable $l_\mu\in \NN_0$. Since $\mu<\gamma$, the order of $\mu$ is less than 
 the order of $\gamma$. Thus, by the induction hypothesis, there exists a polynomial $P_\mu$
with order $<l$ such that 
 $\chi_{\mu}=P_\mu(\chi_{\lambda_1},\ldots,\chi_{\lambda_k})$.
 Similarly, the induction hypothesis gives a decomposition
 $$\chi_{\gamma-\lambda_1}=\chi_{\lambda_1}^{l_1-1}\prod_{i=2}^k\chi_{\lambda_i}^{l_i}+
 P_{\gamma-\lambda_1}(\chi_{\lambda_1},\ldots, \chi_{\lambda_k}),$$
 such that the order of $P_{\gamma-\lambda}$ is smaller than $l-1$.
The result then follows with $P=X_1P_{\gamma-\lambda_1}-\sum_{\mu<\gamma}l_\mu P_\mu$.
\end{proof}

We are now coming back to the special case of the group $G=\Spin(n)$ with $n=2m$ and 
$m\geq 2$.
This group is simply connected and connected and by \cite[Chapter VI, Theorem (6.2)]{BtD}
a system of  \emph{fundamental representations} is given by the representations
$$\Lambda^1,\ldots, \Lambda^{m-2}, \Sigma^+,\Sigma^-$$
defined as follows: 
The representations $\Lambda^i$ act on the complexification $\Lambda^i(\CC^n)$ of the
$i$th exterior power $\Lambda^i(\RR^n)$ by inflating the canonical action of $\SO(n)$ on 
$\Lambda^i(\RR^n)$ to $\Spin(n)$.
Note that these representations extend canonically to $\Pin(n)$ (resp.\  to $\O(n)$, if we view them as 
representations of $\SO(n)$), 
which implies that the $\Lambda^i$ are 
stable (up to equivalence) under conjugation by elements in $\Pin(n)$ (resp. $\O(n)$). It 
is also clear that the $\Lambda^i$ are non-negative, i.e, $\Lambda^i(-1)=1$.

The representations $\Sigma^+, \Sigma^-$ are the half-spin representations on the spaces
$S^+, S^-$ defined as follows: By the isomorphism $\Cl(n)\cong M_{2^m}(\CC)$ 
we find a canonical irreducible action of the complex Clifford algebra $\Cl(n)$ on 
$S:=\CC^{2^m}$. Since $J^2=(-1)^m$, for $J=e_1\cdots e_n$, it follows that 
$\tilde{J}:=i^mJ$ satisfies $\tilde{J}=\tilde{J}^*=\tilde{J}^{-1}$, which implies that
$S$  decomposes into the direct sum of
two orthogonal eigenspaces $S^+, S^-$  for the eigenvalues $\pm 1$ of
 $\tilde{J}$. Since $\tilde{J}x\tilde{J}=JxJ^*=x$ for all $x\in \Cl(n)^0$, these spaces are invariant 
under the action of $\Cl(n)^0$, and then restrict to unitary representations $\Sigma^{\pm}$ of 
$\Spin(n)\subseteq \Cl(n)^0$. One easily checks that conjugation by
$e_1\in \Pin(n)\setminus \Spin(n)$ intertwines these representations. We therefore see 
that $\{\Sigma^+,\Sigma^-\}$ forms one orbit of length two in $\widehat{\Spin(n)}$ under 
conjugation by $\Pin(n)$.  By construction, the representations 
$\Sigma^\pm$ are negative representations, i.e., $\Sigma^\pm(-1)=-1$.

We are now ready to prove the following proposition, which will give the last step in the 
proof of Theorem \ref{prop-general-even} of the previous section.

\begin{proposition}\label{prop-Spin}
Let $n=2m\geq 4$. Then the following are true:
\begin{enumerate}
\item If $x\in \Pin(n)\setminus\Spin(n)$ and $\tau$ is a negative irreducible representation of 
$\Spin(n)$, then $\tau\not\cong x\cdot\tau$. Thus, for the orbit sets $\mathcal O_1$ and $\mathcal O_2$ 
in $\widehat{\Spin(n)}^-$  we get $\mathcal O_1=\emptyset$ and $\mathcal O_2$
is countably infinite.
\item For the action of $\O(n)$ on $\widehat{\SO(n)}$ both orbit sets $\mathcal O_1$ and $\mathcal O_2$ are countably infinite.
\end{enumerate}
\end{proposition}
\begin{proof}
Let  $\chi_1,\ldots,\chi_{m-2}, \chi_+,\chi_-$ denote the characters corresponding to the fundamental
representations $\Lambda^1,\ldots,\Lambda^{m-2},\Sigma^\pm$. 
Let $\tau$  be any negative irreducible representation of $\Spin(n)$ with character $\chi_\tau$.
By Theorem \ref{thm-fund} there exists a unique Polynomial $Q\in \ZZ[X_1,\ldots, X_{m-2}, X_+, X_-]$
such that $\chi_\tau=Q(\chi_1,\ldots,\chi_{m+2},\chi_+,\chi_-)$. By 
 Lemma \ref{lem-polydeco} the polynomial $Q$ is of the form
 $$\left(\prod_{i=1}^{m-2}X_i^{l_i}\right)X_+^{l_+}X_-^{l_-}+P(X_1,\ldots,X_{m-2}, X_+,X_-)$$
 with the order of $P$ less than $l=l_1+\cdots +l_{m-2}+l_++l_-$.
Since $\tau$ is negative, we have $\chi_\tau(-x)=-\chi_\tau(x)$ for all $gx\in \Spin(n)$.
Since for all $x\in \Spin(n)$ we have $\chi_i(-x)=\chi_i(x)$, for all $1\leq i\leq m-2$,
 and $\chi_{\pm}(-x)=-\chi_{\pm}(x)$ we get
$$\chi_\tau=-(-1)^{l_++l_-}\left(\prod_{i=1}^{m-2}\chi_i^{l_i}\right)\chi_+^{l_+}\chi_-^{l_-}+
\tilde{P}(\chi_1,\ldots, \chi_{n-2}, \chi_+,\chi_-)$$
with the order of $\tilde{P}$ less than $l$. By the uniqueness of the polynomial representation 
of $\chi_\tau$ it follows that $(-1)^{l_++l_-}=-1$ and, in particular, that $l_+\neq l_-$.

Suppose now that $x\in \Pin(n)\setminus \Spin(n)$. Since $x\cdot \chi_i=\chi_i$ for all 
$1\leq i\leq m-2$ and $x\chi_+=\chi_-$ (and vice versa) we get
\begin{align*}
x\chi_{\tau}&=Q(x\chi_1,\ldots,x\chi_{m-2},x\chi_+,x\chi_-)\\
&=\left(\prod_{i=1}^{m-2}(x\chi_i)^{l_i}\right)(x\chi_+)^{l_+}(x\chi_-)^{l_-}+
{P}(x\chi_1,\ldots, x\chi_{m-2}, x\chi_+,x\chi_-)\\
&=\left(\prod_{i=1}^{m-2}\chi_i^{l_i}\right)\chi_-^{l_+}\chi_+^{l_-}+
{P}(\chi_1,\ldots, \chi_{m-2}, \chi_-,\chi_+)\\
&=\tilde{Q}(\chi_1,\ldots,\chi_{m-2},\chi_+,\chi_-).
\end{align*}
Since $l_+\neq l_-$ we have $Q\neq \tilde{Q}$, hence $\chi_{x\tau}=x\chi_\tau\neq \chi_{\tau}$,
and therefore $x\tau\not\cong\tau$. This proves (i).

For the proof of (ii) note first that $\widehat{\SO(n)}=\widehat{\Spin(n)}^+$, the set of irreducible  representations $\tau$ of $\Spin(n)$ with $\tau(-1)=1$. Writing its character $\chi_\tau$ 
as $Q(\chi_1,\ldots,\chi_{m+2},\chi_+,\chi_-)$ as above, we see that $\tau\in \widehat{\SO(n)}$ if and 
only if $l_++l_-$ is even. Then a similar computation as above shows that  for $x\in \Pin(n)\setminus\Spin(n)$ we get 
$$x\tau\cong \tau\Leftrightarrow x\chi_\tau=\chi_\tau\Leftrightarrow l_+=l_-.$$
It is  now clear that there are infinitely many representations which are fixed by conjugation and 
there are also infinitely many pairs of conjugate representations in $\widehat{\SO(n)}$.
\end{proof}

\begin{proof}[Proof of Theorem \ref{prop-general-even}]
The proof now follows from the above proposition together with the discussion of 
the even case preceding Proposition \ref{prop-spin}.
\end{proof}

\end{document}